\documentclass[10pt]{amsart}
\usepackage{amsfonts}
\usepackage{amsmath}
\usepackage{amsthm}
\usepackage{amssymb}
\usepackage{latexsym}
\usepackage{multicol}
\usepackage{verbatim}

\advance\textwidth by 1.2in \advance\oddsidemargin by -.6in
\advance\evensidemargin by -.6in
\newtheorem*{cor}{Corollary}
\newtheorem*{lem}{Lemma}
\newtheorem*{prop}{Proposition}

\theoremstyle{definition}
\newtheorem*{defn}{Definition}
\theoremstyle{definition}
\newtheorem*{thm}{Theorem}

\newtheorem*{rem}{Remark}

\newenvironment{pf}{\proof}{\endproof}
\newcounter{cnt}
\newenvironment{enumerit}{\begin{list}{{\hfill\rm(\roman{cnt})\hfill}}{%
\settowidth{\labelwidth}{{\rm(iv)}}\leftmargin=\labelwidth%
\advance\leftmargin by
\labelsep\rightmargin=0pt\usecounter{cnt}}}{\end{list}}

\theoremstyle{remark}

\numberwithin{equation}{section} 

\def\dnl{=\hspace{-.2cm}>\hspace{-.2cm}=}
\def\dnr{=\hspace{-.2cm}<\hspace{-.2cm}=}
\def\sn{-\hspace{-.17cm}-}
\def\wtl{{\rm wt}_\ell}
\def\wt{{\rm wt}}

\def\opl_#1{\text{\scriptsize$\bigoplus\limits_{\text{\footnotesize$#1$}}$}}

\begin{document}

\newcommand{\thmref}[1]{Theorem~\ref{#1}}
\newcommand{\secref}[1]{Section~\ref{#1}}
\newcommand{\lemref}[1]{Lemma~\ref{#1}}
\newcommand{\propref}[1]{Proposition~\ref{#1}}
\newcommand{\corref}[1]{Corollary~\ref{#1}}
\newcommand{\remref}[1]{Remark~\ref{#1}}
\newcommand{\defref}[1]{Definition~\ref{#1}}
\newcommand{\er}[1]{(\ref{#1})}
\newcommand{\id}{\operatorname{id}}
\newcommand{\tensor}{\otimes}
\newcommand{\nc}{\newcommand}
\newcommand{\rnc}{\renewcommand}
\newcommand{\qbinom}[2]{\genfrac[]{0pt}0{#1}{#2}}
\nc{\cal}{\mathcal} \nc{\goth}{\mathfrak} \rnc{\bold}{\mathbf}
\renewcommand{\frak}{\mathfrak}
\newcommand{\desc}{\operatorname{desc}}
\newcommand{\Maj}{\operatorname{Maj}}
\renewcommand{\Bbb}{\mathbb}
\nc\bomega{{\mbox{\boldmath $\omega$}}} \nc\bvpi{{\mbox{\boldmath
$\varpi$}}}
 \nc\balpha{{\mbox{\boldmath $\alpha$}}}

\newcommand{\lie}[1]{\mathfrak{#1}}
\makeatletter
\def\section{\def\@secnumfont{\mdseries}\@startsection{section}{1}%
  \z@{.7\linespacing\@plus\linespacing}{.5\linespacing}%
  {\normalfont\scshape\centering}}
\def\subsection{\def\@secnumfont{\bfseries}\@startsection{subsection}{2}%
  {\parindent}{.5\linespacing\@plus.7\linespacing}{-.5em}%
  {\normalfont\bfseries}}
\makeatother
\def\subl#1{\subsection{}\label{#1}}

\nc{\Cal}{\cal} \nc{\Xp}[1]{X^+(#1)} \nc{\Xm}[1]{X^-(#1)}
\nc{\on}{\operatorname} \nc{\ch}{\mbox{ch}} \nc{\Z}{{\bold Z}}
\nc{\J}{{\cal J}} \nc{\C}{{\bold C}} \nc{\Q}{{\bold Q}}
\renewcommand{\P}{{\cal P}}
\nc{\N}{{\Bbb N}} \nc\boa{\bold a} \nc\bob{\bold b} \nc\boc{\bold
c} \nc\bod{\bold d} \nc\boe{\bold e} \nc\bof{\bold f}
\nc\bog{\bold g} \nc\boh{\bold h} \nc\boi{\bold i} \nc\boj{\bold
j} \nc\bok{\bold k} \nc\bol{\bold l} \nc\bom{\bold m}
\nc\bon{\bold n} \nc\boo{\bold o} \nc\bop{\bold p} \nc\boq{\bold
q} \nc\bor{\bold r} \nc\bos{\bold s} \nc\bou{\bold u}
\nc\bov{\bold v} \nc\bow{\bold w} \nc\boz{\bold z}

\nc\ba{\bold A} \nc\bb{\bold B} \nc\bc{\bold C} \nc\bd{\bold D}
\nc\be{\bold E} \nc\bg{\bold G} \nc\bh{\bold H} \nc\bi{\bold I}
\nc\bj{\bold J} \nc\bk{\bold K} \nc\bl{\bold L} \nc\bm{\bold M}
\nc\bn{\bold N} \nc\bo{\bold O} \nc\bp{\bold P} \nc\bq{\bold Q}
\nc\br{\bold R} \nc\bs{\bold S} \nc\bt{\bold T} \nc\bu{\bold U}
\nc\bv{\bold V} \nc\bw{\bold W} \nc\bz{\bold Z} \nc\bx{\bold X}

\title[Characters and blocks for finite--dimensional representations of
 qaa]{Characters and blocks for finite--dimensional representations of
 quantum affine algebras}

\author{Vyjayanthi Chari and Adriano A. Moura }
\address{Department of Mathematics, University of
California, Riverside, CA 92521.} \email{chari@math.ucr.edu,
adrianoam@math.ucr.edu}

\maketitle

\setcounter{section}{0}

\section*{Introduction}

In this paper we study the category $\cal{C}_q$  of
finite--dimensional representations of a quantum loop algebra
$\bu$. Our aim is to study and to put into a common representation
theoretic framework, two kinds of characters which have been
associated to an object of $\cal{C}_q$. One is the notion of
$q$--characters defined in \cite{FR} which is analogous in this
context, to the usual notion of a character of a
finite--dimensional representation of a simple Lie algebra. The
other, is the notion of  the elliptic character defined in
\cite{EM} which plays the role of the central character for
representations of  semi-simple Lie algebras. Both kinds of
characters are needed in our situation, because the category
$\cal{C}_q$ is not semi-simple and hence the problem of
determining the blocks in this category becomes important.

 The papers \cite{EM} and \cite{FR} use the universal $\cal{R}$--matrix
in fundamental ways to study the elliptic character and the
$q$--character  respectively. In particular, \cite{EM} uses
convergence properties of this matrix and hence, the main result
of the paper describes the blocks in the case  $|q|<1$. Our
methods, which avoid the use of the $\cal R$--matrix,  allows us
to determine the blocks for all $q$ not a root of unity. One of
the conjectures of \cite{FR}, proved in \cite{FM} (see also
\cite{He}) is that the character of simple objects of $\cal C_q$
has a certain cone like form. We prove this result for the quantum
affine algebras associated to a classical Lie algebra in a
representation theoretic way
 rather than  in a combinatorial fashion.  We are actually able to prove a  stronger
version of their result,  which allows us to  give a formula for
the $q$--characters of the fundamental representations in terms of
the braid group action defined in  \cite{C1}.

We now describe the results of this paper.  The algebra $\bu$  has
a large commutative subalgebra $\bu(0)$ and the representations in
$\cal{C}_q$ can be written as a sum of generalized eigenspaces for
this subalgebra. The corresponding eigenvalues are known to be
$n$--tuples of rational functions \cite{FR}, where $n$ is the rank
of the underlying finite--dimensional simple Lie algebra $\lie g$.
We define the $\ell$--weight lattice $\cal{P}_q$ of $\bu$ to be
the multiplicative subgroup  consisting of the invertible elements
in the   ring of $n$--tuples of rational functions in an
indeterminate $u$. It was proved in \cite{C1} that the braid group
of  $\lie g$ acts on $\cal{P}_q$. Using this action, we define in
Section 2 the notion of simple $\ell$--roots and the $\ell$--root
lattice $\cal{Q}_q$. It turns out that $\cal Q_q$ is preserved by
the braid group action. We then give generators and relations for
the quotient group $\Xi_q=\cal{P}_q/\cal{Q}_q$. Our constructions
make sense when $q=1$ and the quotient group $\Xi_1$ is just the
group of functions with finite support  from $\bc^{\times}$ to
$P/Q$, where $P$ and $Q$ are the usual weight and root lattices of
$\lie g$. In Section 7 we prove that the blocks of $\cal{C}_q$ are
in bijective correspondence with elements of $\Xi_q$. A somewhat
unusual feature of  this category, is that  the tensor product of
two blocks is contained in a single block.

Let $\cal P_q^+$ be the monoid in $\cal{P}_q$ generated by
$n$--tuples of polynomials. It was proved in \cite{CPbanff} that
elements of $\cal{P}_q^+$ parametrize the isomorphism classes of
the irreducible objects in $\cal{C}_q$. Motivated by this, we call
the elements of $\cal P_q^+$ the dominant $\ell$-weights. Given
$\bomega\in\cal{P}_q^+$ one can define in a natural way the notion
of an $\ell$--highest weight representation with $\ell$--highest
weight $\bomega$. In \cite{CPweyl}  a family of  universal
$\ell$--highest weight module in $\cal{C}_q$, called the Weyl
module $W(\bomega)$ was constructed and it was  conjectured there
and proved in the case of $sl_2$  that $W(\bomega)$ was isomorphic
to a tensor product of fundamental modules, or in other words,
that $W(\bomega)$ was isomorphic to a standard module. Using some
deep results of Nakajima we deduce this conjecture for a general
simple Lie algebra in section 6.  We then
 prove that  the $\ell$--weights of $W(\bomega)$ and hence,
those of any $\ell$--highest weight representations
 lie in the \lq\lq cone\rq\rq \  $\bomega(\cal{Q}_q^+)^{-1}$,
 here $\cal{Q}_q^+$ is the monoid generated by the simple
 $\ell$--roots. This result  plays an important role in
 Section 7, since it allows us to prove that the Weyl module has
 a well-defined elliptic character.

In sections 4 and 5, we study the $\bu(0)$--decomposition into
generalized eigenspaces of  objects in $\cal{C}_q$. We call these
the $\ell$--weight spaces and the corresponding eigenvalues the
$\ell$--weights of that representation. We prove that the
$\ell$--weights of any object of $\cal C_q$ are in $\cal{P}_q$ and
that the $\ell$--weights determine the usual weights of objects in
$\cal{C}_q$ via a homomorphism $\wt: \cal P_q\to P$. To further
describe the main result of Section 5, it is useful to compare it
with the results  of \cite{FM} and \cite{FR}. In those papers the
authors developed the notion of a $q$--character for objects in
$\cal{C}_q$. In the language of this paper, they are the following
element of $\bz[\cal P_q]$,
$${\rm ch}_\ell(V)=\sum_{\bvpi\in\cal{P}_q}\dim(V_\bvpi)\ e(\bvpi),$$
where $V_\bvpi$ is the generalized
eigenspace of $V$ corresponding to the eigenvalue $\bvpi$, and $\bz[\cal P_q]$ is the integral group ring over $\cal P_q$ with basis elements $e(\bvpi)$.
 The elements
 $A_{i,c}$ of
\cite{FR}, where $i$ varies over the set of simple roots for $\lie
g$ and $c\in\bc^\times$  turn out to correspond to the simple
$\ell$--roots. In \cite{FM} a result corresponding to the main
result of Section 5 of this paper, namely, that the
$\ell$--weights of the irreducible representation $V(\bomega)$
 lie in the cone $\bomega(\cal{Q}_q^+)^{-1}$ was proved using
 combinatorial methods. Our methods  on the other hand are purely
  representation theoretic.
 This allows us
state more precise results on the $\ell$ (or $q$)--character of
the fundamental modules of the classical algebras. Thus, we prove
that the $\ell$--weights of fundamental representations have a
certain invariance under the braid group action analogous to the
invariance of the set of weights under the Weyl group for
finite--dimensional representations of simple Lie algebras.
Interestingly enough, this invariance appears to be a very special
property which fails if the Lie algebra is of exceptional type.
Although there are a number of papers where $q$-characters have
been studied, \cite{He}, \cite{Nak}, there are few explicit
formulas available even  for the fundamental representations,
although there are some conjectures in \cite{FR} and there is a
description of the $\ell$--characters in the $A_n$, $D_n$ case in
terms of tableaux in \cite{Nak}. As an application of our
techniques, we write down the $q$--character of the fundamental
representation corresponding to the adjoint representation when
$\lie g$ is of type $D_n$. The more general case is  studied in
\cite{CM3}.

\vskip12pt
{\bf Acknowledgements:} We thank E. Mukhin and J. Greenstein for useful discussions.

\section{Preliminaries}
In this section, we   recall the definition of quantum affine
algebras and several results on the structure of these algebras.

\subsection{} Let $\frak g$ be a complex finite--dimensional simple
Lie algebra  of rank
 $n$ and let $\frak h$ be a Cartan subalgebra of $\lie g$.
 Set $I=\{1,2,\cdots ,n\}$ and  let $\{\alpha_i:i\in I\}$ (resp.  $\{\omega_i:i\in I\}$)
  be the set of simple roots
 (resp.  fundamental weights) of $\frak g$ with respect to $\frak h$. Let also $\check{\alpha}_i$ denote the simple co-roots. As usual,
 $Q$,
(resp. $P$) denotes the  root (resp. weight) lattice of $\frak g$
and $Q^+$, $P^+$ the non--negative root and weight lattice
respectively. Set $\Gamma=P/Q$.
 Let $A=(a_{ij})_{i,j\in I}$  be the $n\times n$ Cartan matrix of $\frak g$ and let  $\check h\in\bz$ be the
 dual Coxeter number of $\lie g$.  Fix non--negative integers $d_i$ such that
the matrix $(d_ia_{ij})$ is symmetric. Assume that the nodes of
the Dynkin diagram of $\lie g$ are numbered as in Table 1 below and let $I_\bullet$
denote the subset of $I$ consisting of the shaded nodes.

\newpage

{\centerline {\bf Table 1}}
\begin{multicols}{2}
\begin{itemize}

\item $A_n$ : {\large \vspace{-.55cm}
$$\stackrel{1}{\bullet}\hspace{-.18cm}\sn\hspace{-.18cm}\stackrel{2}{\circ}
\dots
\stackrel{\text{n-1}}{\circ}\hspace{-.29cm}\sn\hspace{-.18cm}\stackrel{\text{n}}{\circ}$$}

\item $B_n$ : {\large \vspace{-.55cm}
$$\stackrel{1}{\circ}\hspace{-.18cm}\sn\hspace{-.18cm}\stackrel{2}{\circ}
\dots
\stackrel{\text{n-1}}{\circ}\hspace{-.18cm}\dnl\hspace{-.07cm}\stackrel{\text{n}}{\bullet}$$}

\item $C_n$ : {\large \vspace{-.55cm}
$$\stackrel{1}{\bullet}\hspace{-.18cm}\sn\hspace{-.18cm}\stackrel{2}{\circ}
\dots
\stackrel{\text{n-1}}{\circ}\hspace{-.18cm}\dnr\hspace{-.07cm}\stackrel{\text{n}}{\circ}$$}

\item $D_n$, $n$ odd : {\large \vspace{-.75cm}
$$\hspace{.95cm}\stackrel{1}{\circ}\hspace{-.18cm}\sn\hspace{-.18cm}\stackrel{2}{\circ}
\dots
\stackrel{\text{n-2}}{\circ}\hspace{-.29cm}\sn\hspace{-.28cm}\stackrel{\text{n-1}}{\circ}$$
\vspace{-.82cm}$$\hspace{1.8cm}|$$
\vspace{-.73cm}$$\hspace{2.11cm}\bullet\text{ \small n}$$}

\item $D_n$, $n$ even : {\large \vspace{-.75cm}
$$\hspace{.95cm}\stackrel{1}{\circ}\hspace{-.18cm}\sn\hspace{-.18cm}\stackrel{2}{\circ}
\dots
\stackrel{\text{n-2}}{\circ}\hspace{-.29cm}\sn\hspace{-.28cm}\stackrel{\text{n-1}}{\bullet}$$
\vspace{-.82cm}$$\hspace{1.8cm}|$$
\vspace{-.73cm}$$\hspace{2.11cm}\bullet\text{ \small n}$$}
\\
\item $E_6$  : {\large \vspace{-.55cm}
$$\stackrel{1}{\bullet}\hspace{-.18cm}\sn\hspace{-.18cm}\stackrel{2}{\circ}
\hspace{-.18cm}\sn\hspace{-.18cm}
\stackrel{\text{3}}{\circ}\hspace{-.18cm}\sn\hspace{-.18cm}\stackrel{4}{\circ}\hspace{-.18cm}\sn
\hspace{-.18cm}\stackrel{5}{\circ}$$
\vspace{-.75cm}$$\hspace{-.00cm}|$$
\vspace{-.75cm}$$\hspace{.28cm}\circ\text{ \footnotesize$6$}$$}

\item $E_7$  : {\large \vspace{-.65cm}
$$\hspace{.5cm}\stackrel{1}{\bullet}\hspace{-.18cm}\sn\hspace{-.18cm}\stackrel{2}{\circ}
\hspace{-.18cm}\sn\hspace{-.18cm}
\stackrel{\text{3}}{\circ}\hspace{-.18cm}\sn\hspace{-.18cm}
\stackrel{4}{\circ}\hspace{-.18cm}\sn
\hspace{-.18cm}\stackrel{5}{\circ}\hspace{-.18cm}
\sn\hspace{-.18cm}\stackrel{6}{\circ}$$
\vspace{-.82cm}$$\hspace{.98cm}|$$
\vspace{-.75cm}$$\hspace{1.26cm}\circ\text{ \footnotesize$7$}$$}

\item $E_8$  : {\large \vspace{-.75cm}
$$\hspace{1cm}\stackrel{1}{\bullet}\hspace{-.18cm}\sn\hspace{-.18cm}\stackrel{2}{\circ}
\hspace{-.18cm}\sn\hspace{-.18cm}
\stackrel{\text{3}}{\circ}\hspace{-.18cm}
\sn\hspace{-.18cm}\stackrel{4}{\circ}\hspace{-.18cm}\sn
\hspace{-.18cm}\stackrel{5}{\circ}\hspace{-.18cm}
\sn\hspace{-.18cm}\stackrel{6}{\circ}
\hspace{-.18cm}\sn\hspace{-.18cm}\stackrel{7}{\circ}$$
\vspace{-.82cm}$$\hspace{1.95cm}|$$
\vspace{-.75cm}$$\hspace{2.23cm}\circ\text{ \footnotesize$8$}$$}

\item $F_4$ : {\large \vspace{-.55cm}
$$\hspace{-.4cm}\stackrel{1}{\bullet}\hspace{-.18cm}\sn\hspace{-.18cm}\stackrel{2}{\circ}
\hspace{-.08cm}\dnr\hspace{-.08cm}\stackrel{3}{\circ}
\hspace{-.18cm}\sn\hspace{-.18cm}\stackrel{4}{\circ}$$}

\item $G_2$ : {\large \vspace{-.55cm}
$$\hspace{-1.4cm}\stackrel{1}{\bullet}\hspace{-.12cm}\equiv\hspace{-.18cm}<\hspace{-.18cm}\equiv
\hspace{-.1cm}\stackrel{2}{\circ}$$}

\end{itemize}

\end{multicols}

\subl{} Denote by  $W$  the Weyl group of $\lie g$, then $W$ is
generated by simple reflections $\{s_i:i\in I\}$. For $w\in W$,
let $\ell(w)$ denote the length of a reduced expression for $w$.
Let $w_0$ denote the longest element of $W$, then $w_0$ defines a
permutation of $I$, given by $w_0\alpha_i=-\alpha_{w_0i}$. Given
$\lambda=\sum_{i\in I} \lambda_i\omega_i\in P^+$, let
$I(\lambda)=\{i\in I:\lambda_i=0\}$ and let $W(\lambda)$ be the
subgroup of $W$ generated by $\{s_i: i\in I(\lambda)\}$. The
following lemma is well--known \cite{H2}.
\begin{lem}\label{wlambda} Let $\lambda=\sum_{i\in I} \lambda_i\omega_i\in P^+$. Then,
\begin{enumerit}
\item[(i)] $W(\lambda)=\{w\in W: w\lambda=\lambda\}$.
\item[(ii)] Each left (right) coset of $W(\lambda)$ in $W$ contains a unique element
of minimal length. Denote by $W_\lambda$ the set of left coset
representatives of minimal length.
\item[(iii)] Suppose that $w\in W_\lambda$ and that $w=s_jw'$ for
some $w'\in W$ with $\ell(w')=\ell(w)-1$. Then, $w'\in W_\lambda$.
\end{enumerit} \hfill\qedsymbol
\end{lem}

\subl{}
 The next lemma is easily checked using the explicit
formulas for the fundamental weights given in \cite{H}.
\begin{lem}\label{domfund} Suppose that $\lie g$ is of type $A_n$, $B_n$, $C_n$
or $D_n$.
\begin{enumerit}
\item[(i)] If  $\lambda\in P^+$ is such that
$\omega_i-\lambda\in Q^+$ for some $i\in I$, then either $\lambda
=0$ or $\lambda=\omega_j$ for some $j\in I$ with $j\leq i$.
\item[(ii)] Let $i,j\in I$ and assume that $i>j$. Then,
$$\omega_i-\omega_j-\alpha_j\notin Q^+,\ \
\omega_i-\omega_j-2\alpha_{j+1}\notin Q^+.$$
\end{enumerit}\hfill\qedsymbol\end{lem}

\subsection{} Let $q\in\bc^\times$ and assume that $q$ is not a root
of unity. For $r,m\in\bn$, $m\ge r$, define complex numbers,
\begin{equation*}
[m]_q=\frac{q^m -q^{-m}}{q -q^{-1}},\ \ \ \ [m]_q!
=[m]_q[m-1]_q\ldots [2]_q[1]_q,\ \ \ \ \left[\begin{matrix} m\\
r\end{matrix}\right]_q = \frac{[m]_q!}{[r]_q![m-r]_q!}.
\end{equation*}
 Set $q_i=q^{d_i}$ and $[m]_i=[m]_{q_i}$.
The quantum loop algebra  $\bu$ of $\lie g$ is  the algebra with
generators $x_{i,r}^{{}\pm{}}$ ($i\in I$, $r\in\bz$), $K_i^{{}\pm
1}$ ($i\in I$), $h_{i,r}$ ($i\in I$, $r\in \bz\backslash\{0\}$)
and the following defining relations:
\begin{align*}
   K_iK_i^{-1} = K_i^{-1}K_i& =1, \ \
 K_iK_j =K_jK_i,\\  K_ih_{j,r}& =h_{j,r}K_i,\\
 K_ix_{j,r}^\pm K_i^{-1} &= q_i^{{}\pm
    a_{ij}}x_{j,r}^{{}\pm{}},\ \ \\
  [h_{i,r},h_{j,s}]=0,\; \; & [h_{i,r} , x_{j,s}^{{}\pm{}}] =
  \pm\frac1r[ra_{ij}]_{q^i}x_{j,r+s}^{{}\pm{}},\\
 x_{i,r+1}^{{}\pm{}}x_{j,s}^{{}\pm{}} -q_i^{{}\pm
    a_{ij}}x_{j,s}^{{}\pm{}}x_{i,r+1}^{{}\pm{}} &=q_i^{{}\pm
    a_{ij}}x_{i,r}^{{}\pm{}}x_{j,s+1}^{{}\pm{}}
  -x_{j,s+1}^{{}\pm{}}x_{i,r}^{{}\pm{}},\\ [x_{i,r}^+ ,
  x_{j,s}^-]=\delta_{i,j} & \frac{ \psi_{i,r+s}^+ -
    \psi_{i,r+s}^-}{q_i - q_i^{-1}},\\
\sum_{\pi\in\Sigma_m}\sum_{k=0}^m(-1)^k\left[\begin{matrix}m\\k\end{matrix}
\right]_{i}
  x_{i, r_{\pi(1)}}^{{}\pm{}}\ldots x_{i,r_{\pi(k)}}^{{}\pm{}} &
  x_{j,s}^{{}\pm{}} x_{i, r_{\pi(k+1)}}^{{}\pm{}}\ldots
  x_{i,r_{\pi(m)}}^{{}\pm{}} =0,\ \ \text{if $i\ne j$},
\end{align*}
for all sequences of integers $r_1,\ldots, r_m$, where $m
=1-a_{ij}$, $\Sigma_m$ is the symmetric group on $m$ letters, and
the $\psi_{i,r}^{{}\pm{}}$ are determined by equating powers of
$u$ in the formal power series $$\sum_{r=0}^{\infty}\psi_{i,\pm
r}^{{}\pm{}}u^{{}\pm r} = K_i^{{}\pm 1}
{\text{exp}}\left(\pm(q_i-q_i^{-1})\sum_{s=1}^{\infty}h_{i,\pm s}
u^{{}\pm s}\right).$$ \vskip 12pt

\subsection{}
For $i\in I$,  set
\begin{equation*}
h^\pm_i(u)=\sum_{k=1}^\infty \frac{q^{\pm k}h_{i,\pm
k}}{[k]_{i}}u^k,\end{equation*} and define elements $P_{i,\pm k}$,
$i\in I$, $k\in\bz$, $k\ge 0$, by the generating series,
\begin{equation}\label{hp}
P^\pm_i(u)=\sum_{k=0}^\infty P_{i,\pm
k}u^k=\exp(-h^\pm_i(u)).\end{equation}

Let $\bu^\pm(0)$ be the subalgebra of $\bu$ generated by the
elements $h_{i,\pm k}$ $i\in I$, $k\in\bz$, $k>0$. It is easy to
see that $\bu^\pm(0)$ are  commutative subalgebras of $\bu$ and
that monomials in the $h_{i,\pm k}$, $i\in I$, $k\in\bz$ form a
basis of $\bu^\pm(0)$. Notice also that $\bu^\pm(0)$ is also the
subalgebra generated by the elements $P_{i,\pm k}$, $i\in I$,
$k\in\bz$, $k>0$.  The subalgebra of $\bu$ generated by
$x_{i,0}^\pm$ and $K_i^{\pm 1}$, $i\in I$ is isomorphic to the
quantized enveloping algebra $\bu^{fin}$ of $\lie g$. For each
$i\in I$,  the subalgebra $\bu^i$ of $\bu$ generated by the
elements $x^\pm_{i,r}$, $K_i^{\pm 1}$, $h_{i,s}$, $r\in\bz$, $0\ne
s\in\bz$ is isomorphic to the quantum loop algebra  of $sl_2$.

\subsection{} It is well--known that $\bu$ has the structure of
a Hopf algebra,  let $\Delta:\bu\to\bu\otimes \bu$ be the
comultiplication. Although, explicit formulas for the
comultiplication on the generators $x_{i,k}^\pm$, $h_{i,r}$, are
not known, the next proposition proved in \cite[Proposition
5.3]{B}, \cite{Da}, \cite[Proposition 5.4]{BCP} gives partial
information that suffices for our purposes.

Define subspaces $X^\pm$ of $\bu$, by
\begin{equation*}
X^\pm =\sum_{i\in I, r\in\bz}  \bc x_{i,r}^\pm.
\end{equation*}

\begin{prop}\label{comultip}
 Modulo $\bu X^-\otimes \bu X^+$, we have
\begin{eqnarray*}&
\Delta(h_{i,s})& =h_{i,s}\otimes 1+1\otimes h_{i,s} \ \
(s\in\bz_+, s>0),\\ &\Delta(P_{i,r}) &=\sum_{m=0}^r
P_{i,r-m}\otimes P_{i,m},\ \ (r\in\bz_+, r>0 ).
\end{eqnarray*}

\end{prop}

\section{ Braid group actions, the  $\ell$--weight lattice and the  $\ell$--root lattice} In this
section we introduce the notion of $\ell$--integral weights and
$\ell$--roots. These are certain multiplicative subgroups of the
ring of algebra homomorphisms  $\text{Hom}(\bu(0),\bc)$.

\subl{ell1} Set $$\cal {A}=\{f\in\bc[[u]]: f(0) =1\},$$ where
$\bc[[u]]$ is the ring of formal power series in an indeterminate
$u$.  Clearly $\cal{A}$ is a group under multiplication. Given
$f\in\cal{A}$ and $r\in\bz^+$,  we let $f_r$ denote the
coefficient of $u^r$ in $f$. Let $$\ba=\cal{A}\cap\bc(u).$$ Then
$\ba$ is a free subgroup of $\cal{A}$ with generators
$\{1-au:a\in\bc^\times\}$. Given $f^+\in\bc[u]$ with $f^+(0)=1$,
define $f^-\in\bc[u]$ by $$f^-(u) = u^{\text{deg} f^+}f^+(u^{-1})/
(u^{\text{deg} f^+}f^+(u^{-1}))|_{u=0}.$$ Given
$\varpi^+=f^+/g^+\in\ba$ set $\varpi^- =f^-/g^-$.

Define an  injective group homomorphism $\iota: \ba^n\to
\text{Hom}(\bu(0),\bc)$ by extending, $$ \iota(\bvpi)(P_i^\pm(u))=
\varpi_i^\pm ,$$ where $\bvpi=(\varpi^+_1,\cdots
,\varpi^+_n)\in\ba^n$ and the equality is one of power series. We
call $\iota(\ba^n)$ the $\ell$--integral weight lattice of $\bu$
and henceforth denote it by $\cal{P}_q$.  In what follows we shall
identify $\cal{P}_q$ and $\ba^n$ and denote elements of
$\cal{P}_q$ as $n$--tuples of elements from $\ba$. For $i\in I$
and $a\in\bc^\times$, let $\bomega_{i,a}\in\cal{P}_q$ denote the
element whose  $i^{th}$ entry is $1-au$ and all other entries 1.
We call these the $\ell$--fundamental weights. It is obvious that
$\cal{P}_q$ is generated freely as an abelian group by the
elements $\bomega_{i,a}$, $i\in I$, $a\in\bc^\times$. Let
$\cal{P}_q^+$ denote the monoid generated by $1$ and the elements
$\bomega_{i,a}$, $i\in I$, $a\in\bc^\times$, clearly $\cal{P}_q^+$
is isomorphic to the monoid in $\ba^n$ consisting of $n$--tuples
of polynomials with constant term one and we call such elements
$\ell$--dominant weights.
\begin{defn} Let $\wt:\cal{P}_q\to P$ be the group
homomorphism defined by extending,
$$\wt(\bomega_{i,a})=\omega_i.$$
\end{defn}

\subl{ell2} We now define an action of a  braid group on
$\cal{P}_q$. Let $\cal{B}$ be the group generated by elements
$T_i$ ($i\in I$) with defining relations:
\begin{align*}
T_iT_j &=T_jT_i,\ \ \text{if}\ \ a_{ij} =0,\\
T_iT_jT_i& =T_jT_iT_j,\ \ \text{if}\ \ a_{ij}a_{ji} =1,\\
(T_iT_j)^2&= (T_jT_i)^2,\ \ \text{if}\ \ a_{ij}a_{ji}=2,\\
(T_iT_j)^3&= (T_jT_i)^3,\ \ \text{if}\  \ a_{ij}a_{ji}
=3,\end{align*} where $i,j\in\{1,2,\cdots ,n\}$.
 The next proposition is a reformulation of
 \cite[Proposition 3.1]{C1} and can be easily checked.
\begin{prop}\label{braid} The following formulas define an action of
$\cal{B}$ on $ \cal{P}_q$: let $\bvpi=(\varpi_1,\cdots ,\varpi_n)
\in\cal{P}_q$, then $T_i(\bvpi)$ is defined by,
\begin{align*}
(T_i\bvpi)_j & =\varpi_j,\ \  {\text{if}}\ a_{ji}=0,\\
(T_i\bvpi)_j & =\varpi_j(u)\varpi_i(q_iu),\ \  {\text{if}}\
a_{ji}=-1,\\ (T_i\bvpi)_j & =\varpi_j(u)\varpi_i(q^3u)
\varpi_i(qu),\ \ {\text{if}}\ a_{ji}=-2,\\ (T_i\bvpi)_j &=
\varpi_j(u)\varpi_i(q^5u)\varpi_i(q^3u)\varpi_i(qu),\ \
{\text{if}}\ a_{ji}=-3,\\
 (T_i\bvpi)_i&
 =\frac{1}{\varpi_i(q_i^2u)}.\end{align*}\hfill\qedsymbol
\end{prop}

 \subl{}For $i\in I$, set $$\balpha_{i,a}=
(T_{i}(\bomega_{i, a}))^{-1}\bomega_{i, a}. $$ Clearly we have
$\balpha_{i,a}\in \cal{P}_q$ and we let $\cal{Q}_q$ be the
subgroup of $\cal{P}_q$ generated by the $\balpha_{i,a}$. We call
$\balpha_{i,a}$ the $\ell$--simple roots and $\cal{Q}_q$ the
$\ell$--root lattice and let $\cal{Q}_q^+$ be  the monoid
generated by $\balpha_{i,a}$, $i\in I$, $a\in\bc^\times$, and $\cal Q_q^-=(\cal Q_q^+)^{-1}$.
 Given $f\in\bc[u]$, say
$f=(1-a_1u)\cdots (1-a_ru)$, $a_1,\cdots ,a_r\in\bc^\times$, set
$$\balpha_{i,f}=\prod_{m=1}^r\balpha_{i,a_m},$$ and
$$\balpha_{i,f/g}=\balpha_{i,f}(\balpha_{i,g})^{-1}.$$ Finally,
given $\bvpi\in\ba^n$, set
$$\balpha_{i,\bvpi}=\balpha_{i,(\bvpi)_i}.$$ It is now clear that
Proposition \ref{braid} is equivalent to,
\begin{equation}\label{braid2}
T_i(\bvpi)=\bvpi(\balpha_{i,\bvpi})^{-1},\ \ \forall\
\bvpi\in\ba^n.
\end{equation}
In particular, we have
\begin{prop}
The action of the braid group on $\cal{P}_q$ preserves
 $\cal{Q}_q$.
\hfill\qedsymbol
\end{prop}

\subl{} We list the simple roots for the various classical Lie
algebras below for the reader's convenience.

If $\lie g$ is of type $A_n$, then
$$\balpha_{i,a}(u)=\bomega_{i-1,aq}^{-1}\bomega_{i,a}\bomega_{i,aq^2}\bomega_{i+1,
aq}^{-1}, \ \ i\in I$$ where we understand
$\bomega_{-1,a}=\bomega_{n+1,a}=1$.

If $\lie g$ is of type $B_n$, then
\begin{eqnarray*}&\balpha_{i,a}=(\bomega_{i-1,aq^2})^{-1}\bomega_{i,a}\bomega_{i,aq^4}(\bomega_{i+1,
aq^2})^{-1},\ \  i\in I, \ i\ne n-1,n,\\
 &\balpha_{n-1,a}=
(\bomega_{n-2,aq^2})^{-1}\bomega_{n-1,a}\bomega_{n-1,aq^4}(\bomega_{n,aq}\bomega_{n,aq^3})^{-1},\\
&\balpha_{n,a}
=(\bomega_{n-1,aq})^{-1}\bomega_{n,a}\bomega_{n,aq^2}.\end{eqnarray*}

If $\lie g$ is of type $C_n$, then
\begin{eqnarray*}&\balpha_{i,a}=(\bomega_{i-1,aq})^{-1}\bomega_{i,a}\bomega_{i,aq^2}(\bomega_{i+1,
aq})^{-1},\ \  i\in I, \ i\ne n,\\ &\balpha_{n,a}
=(\bomega_{n-1,aq}\bomega_{n-1aq^3})^{-1}\bomega_{n,a}\bomega_{n,aq^4}.\end{eqnarray*}

If $\lie g$ is of type $D_n$, then
\begin{eqnarray*}&\balpha_{i,a}=(\bomega_{i-1,aq})^{-1}\bomega_{i,a}\bomega_{i,aq^2}(\bomega_{i+1,
aq})^{-1},\ \  i\in I, \ i\ne n-2, n-1,n,\\
 &\balpha_{n-2,a}=
(\bomega_{n-3,aq})^{-1}\bomega_{n-2,a}\bomega_{n-2,aq^2}(\bomega_{n-1,aq}\bomega_{n,aq})^{-1},\\
&\balpha_{n-1,a}
=(\bomega_{n-2,aq})^{-1}\bomega_{n-1,a}\bomega_{n-1,aq^2},\\
&\balpha_{n,a}
=(\bomega_{n-2,aq})^{-1}\bomega_{n,a}\bomega_{n,aq^2}
.\end{eqnarray*}

\begin{rem} The elements $\balpha_{i,a}$ are essentially the
elements $A_{i,a}$ defined in \cite{FR}.\end{rem}

\subl{}  Let $w\in W$ and assume that $w=s_{i_1}s_{i_2}\cdots
s_{i_k}$ is a reduced expression for $w$.
 Set $T_w=T_{i_1}\cdots T_{i_k}$.
 It is well--known that $T_w$ is independent of the choice of the reduced expression for $w$. Given
$\bvpi\in\ba^n$ and $w\in W$, we have
\begin{equation*} T_w\bvpi=T_{i_1}T_{i_2}\cdots T_{i_k}\bvpi =\left((T
_w\bvpi)_1, \cdots ,(T_w\bvpi)_n\right).\end{equation*}
\begin{lem}\label{wconjt}\hfill
 \begin{enumerit}
\item[(i)] For all $w\in W$, and $\bvpi\in \cal{P}_q$, we have
$$\wt(T_w\bvpi)=w(\wt(\bvpi)).$$
\item[(ii)]
Suppose that $\bvpi_r\in\cal{P}_q$, $r=1,2$ are such that
$\wt(\bvpi_r)\in P^+$, $r=1,2$. Then
$$T_{w_1}\bvpi_1=T_{w_2}\bvpi_2\implies
\wt(\bvpi_1)=\wt(\bvpi_2).$$ Further if we set
$\lambda=\wt(\bvpi_1)$ then $$w_1^{-1}w_2\in W(\lambda).$$
\end{enumerit}
\end{lem}
\begin{pf} It suffices to check the result when $T_w=T_j$ and
$\bvpi=\bomega_{i,a}$ for some $i,j\in I$ and $a\in\bc^\times$.
But this is now immediate from Proposition \ref{braid}.  To prove
(ii), notice that (i) implies that  $$w_1(\wt(\bvpi_1))=
w_2(\wt(\bvpi_2)).$$ Since $\wt(\bvpi_r)\in P^+$, this implies
that $\wt(\bvpi_1)=\wt(\bvpi_2)$ and hence also that
$w_1^{-1}w_2\in W(\lambda).$
\end{pf}

\subl{}

The following proposition can be deduced from the results of
\cite{C1}, using some representation theory. However, an
elementary case by case proof is sketched below.
\begin{prop}\label{longest}
 Let $\bomega\in\cal{P}_q^+$. Then,
$(T_{w_0}\bomega)^{-1}\in\cal{P}_q^+$. More precisely, we have
$$(T_{w_0}\bomega)_i=((\bomega)_{w_0i}(q^{\check h}u))^{-1}.$$
\end{prop}
\begin{pf} Assume first that $\lie g$ is of type $A_n$,
then $w_0$ has a reduced expression of the form$$w_0=s_1\cdots
s_ns_1\cdots s_{n-1}\cdots s_1s_2s_1.$$ Proceed by induction on
$n$, noting that induction obviously begins at $n=1$. Since
$$w_0'=s_1\cdots s_{n-1}\cdots s_1s_2s_1$$
 is the reduced expression
for the longest element of $A_{n-1}$, we can assume by induction
that $$T_{w_0'}(\bomega)_j=(\bomega)_{n-j}(q^nu)^{-1},\ \ 1\le
j\le n-1, \
T_{w_0'}(\bomega)_n=(\bomega)_{1}(q^{n-1}u)\bomega_n.$$ A simple
computation now gives the result. The proof for the other Lie
algebras is similar working with the reduced expressions for $w_0$
given in \cite[Section 6]{C1}.
\end{pf}

\subl{}
\begin{lem}\label{wconj} Suppose that $\bvpi,\bvpi'\in\cal{P}_q$ are such that
$\bvpi'=T_w\bvpi$ for some $w\in W$. Then $$\bvpi(\bvpi')^{-1}\in
\cal{Q}_q.$$ Moreover if $\bvpi\in\cal{P}_q^\pm$, then
$$\bvpi(\bvpi')^{-1}\in \cal{Q}^\pm_q.$$
\end{lem}
\begin{pf} It suffices to prove the lemma when $w=s_i$ and
$\bvpi=\bomega_{j,a}$for some $i,j\in I$. If $j=i$, the result is
immediate from the definition of $\balpha_{i,a}(u)$. If $j\ne i$,
then the result is immediate since
$T_i(\bomega_{j,a})=\bomega_{j,a}$.
\end{pf}

\subl{} \begin{lem}\label{coset} Given $\bvpi\in\cal{P}_q$ there
exists $\bomega\in\cal{P}_q^+$ such that
$\bomega(\bvpi)^{-1}\in\cal{Q}_q^+$.
\end{lem}
\begin{pf} It clearly suffices to prove the lemma in the case when
$\bvpi=\bomega_{i,a}^{-1}$ for some $i\in I$, $a\in\bc^\times$. By
Proposition \ref{longest} we see that
$T_{w_0}(\bomega_{i,a})^{-1}=\bomega_{w_0i,aq^{\check h}}$. The result is
now immediate from Lemma \ref{wconj}.
\end{pf}

\section{The group $\cal{P}_q/\cal{Q}_q$}
 In this section we give  a set of generators and relations
for the quotient group $\cal{P}_q/\cal{Q}_q$ in the case when
$\lie g$ is of type $A_n$, $B_n$, $C_n$ or $D_n$. The case of the
exceptional algebras is postponed to the appendix.

\subl{}
\begin{prop}\label{gammaq} \hfill
\begin{enumerit}
\item[(i)] Assume that $\lie g$ is of type $A_n$, $B_n$, $C_n$ or
$D_m$, where $m$ is odd.
 The group $\cal{P}_q/\cal{Q}_q$ is isomorphic to the
(additive) abelian group $\Xi_q$ with generators
$\{\chi_a:a\in\bc^\times\}$ and relations:
\begin{alignat*}{3}
& \sum_{r=0}^{n}\chi_{aq^{n+1-2r}} =0, & &\text{ if } \lie g =
A_n,\\ & \chi_a+\chi_{aq^{4n-2}} =0, & &\text{ if } \lie g =
B_n,\\ & \chi_a+\chi_{aq^{2n+2}} =0, & &\text{ if } \lie g =
C_n,\\ & \chi_a+\chi_{aq^{2}}+\chi_{aq^{2m-2}}+\chi_{aq^{2m}}=0, &
&\text{ if } \lie g = D_m,
\end{alignat*}
for all $a\in\bc^\times$.
\item[(ii)] If $\lie g$ is of type $D_m$ with $m$ even, then $\cal{P}_q/\cal{Q}_q$ is isomorphic to the
(additive) abelian group with generators
$\{\chi^\pm_a:a\in\bc^\times\}$ and relations:
\begin{equation*} \chi^\pm_a+\chi^\pm_{aq^{2m-2}}=0,\ \
\chi^-_a+\chi^-_{aq^2}+\chi^+_{aq^{2m-2}}+\chi^+_{aq^{2m}}=0,\end{equation*}
for all $a\in\bc^\times$.
\end{enumerit}
\end{prop}

The proposition is proved in the rest of the section.

\subl{}
 We begin with
the following Lemma whose proof is an obvious computation.
\begin{lem} \label{ellfund}\hfill
\begin{enumerit} \item[(i)] Suppose that $\lie g$ is of type $A_n$
or $C_n$ (resp. $B_n,D_n$). Then, for all $i\in I$, (resp $i\ne
n$, $i\ne n-1,n$) we have,
\begin{eqnarray*}&\bomega_{1,aq_1^{i-1}}\bomega_{i-1,aq_{i-1}^{-1}}=\bomega_{i,a}
\left(\prod_{j=1}^{i-1}\balpha_{j,aq_j^{i-j-2}}\right)
.\end{eqnarray*} Further,
$$\bomega_{1,aq_1^{-n-1}}\bomega_{n,a}=\prod_{j=1}^n\balpha_{j,
aq_j^{-n+j-2}},\ \ \text{if}\  \lie g =A_n,$$ and
$$\bomega_{1,aq_1^{n+1}}\bomega_{1,aq_1^{-n-1}}=\left(\prod_{j=1}^{n-1}\balpha_{j,aq_j^{n-j}}
\balpha_{j,aq_j^{-n+j-2}}\right)\balpha_{n,aq_j^{-2}},\ \
\text{if} \ \lie g=C_n.$$
\item[(ii)] Assume that $\lie g$ is of type $B_n$. Then,
$$\bomega_{n,aq_n^{2i-1}}\bomega_{n,aq_n^{-2i+1}}=\bomega_{n-i,a}\left(\prod_{j=n-i+1}^{n-1}
\prod_{r=0}^{j-(n-i+1)} \balpha_{j,aq_n^{2(n-i-j)+4r}}\right)
\left(\prod_{r=0}^{i-1}\balpha_{n,aq_n^{-2i+1+4r}}\right).$$
Further,
$$\bomega_{n,aq_n^{2n-1}}\bomega_{n,aq_n^{-(2n-1)}}=\left(\prod_{j=1}^{n-1}
\prod_{r=0}^{j-1} \balpha_{j,aq_j^{-2j+4r}}\right)
\left(\prod_{r=0}^{n-1}\balpha_{n,aq_n^{-2n+1+4r}}\right).$$
\item[(iii)] Assume that $\lie g$ is of type $D_n$. Then, for $j=1,2,\cdots, [(n-1)/2]$,
\begin{align*}
 \bomega_{n,aq^{2j-1}}\bomega_{n,aq^{-2j+1}}= &\ \bomega_{n-2j,a}
\left(\prod_{k=n-2j+1}^{n-2}\prod_{r=0}^{k-(n-2j+1)}
\balpha_{k,aq^{n-2j-k+2r}}\right)\times\\
&\left(\prod_{r=0}^{j-2}\balpha_{n-1,aq^{-2j+3+4r}}\right)
 \left(\prod_{r=0}^{j-1}\balpha_{n,aq^{-2j+1+4r}}\right).
\end{align*}
Also, for $j=1,2,\cdots, [(n-2)/2]$,
{\samepage \begin{align*}
 \bomega_{n-1,aq^{2j}}\bomega_{n,aq^{-2j}}=\bomega_{n-2j-1,a}
\left(\prod_{k=n-2j}^{n-2}\prod_{r=0}^{k-(n-2j)}
\balpha_{k,aq^{n-2j-1-k+2r}}\right)
\left(\prod_{r=0}^{j-1}\balpha_{n-1,aq^{-2j+2+4r}}\balpha_{n,aq^{-2j+4r}}\right).
\end{align*}
Similar formulas hold interchanging
$n$ and $n-1$ on the left hand side.} In addition, if  $n$ is odd,
we have
\begin{align*}
\bomega_{n,aq^2}\bomega_{n,a}\bomega_{n,aq^{-2n+4}}=& \ \bomega_{n-1,a}
\left(\prod_{k=1}^{n-2}\prod_{r=0}^{k-1}
\balpha_{k,aq^{3-n-k+2r}}\right)
\left(\prod_{r=0}^{\frac{n-5}{2}}\balpha_{n-1,aq^{6-2n+4r}}\right)\times\\
&\left(\balpha_{n,a}\prod_{r=0}^{\frac{n-3}{2}}\balpha_{n,aq^{4-2n+4r}}\right),
\end{align*}
\begin{align*}
\bomega_{n-1,aq^{n-1}}\bomega_{n,aq^{-(n-1)}} = &\
\left(\prod_{k=1}^{n-2}\prod_{r=0}^{k-1}
\balpha_{k,aq^{-k+2r}}\right)
\left(\prod_{r=0}^{\frac{n-3}{2}}\balpha_{n-1,aq^{3-n+4r}}\balpha_{n,aq^{1-n+4r}}\right),
\end{align*}
\begin{align*}
\bomega_{n,a}\bomega_{n,aq^{2}}\bomega_{n,aq^{2n-2}}\bomega_{n,aq^{2n}}
= &\ \left(\prod_{k=1}^{n-2}\prod_{r=0}^{k-1}
\balpha_{k,aq^{n-k-1+2r}} \balpha_{k,aq^{n+k-1-2r}}\right)
\left(\prod_{r=0}^{n-3}\balpha_{n-1,aq^{2+2r}}\right)\times\\
&\left(\prod_{r=0}^{n-1}\balpha_{n,aq^{2r}}\right).
\end{align*}
 If $n$ is  even we have
\begin{align*}
\bomega_{n,aq^{n-1}}\bomega_{n,aq^{-(n-1)}} = &\
\left(\prod_{k=1}^{n-2}\prod_{r=0}^{k-1}
\balpha_{k,aq^{-k+2r}}\right)
\left(\prod_{r=0}^{\frac{n-4}{2}}\balpha_{n-1,aq^{3-n+4r}}\right)
\left(\prod_{r=0}^{\frac{n-2}{2}}\balpha_{n,aq^{1-n+4r}}\right)
\end{align*}
and
\begin{align*}
\bomega_{n-1,a}\bomega_{n-1,aq^{2}}\bomega_{n,aq^{2n-2}}\bomega_{n,aq^{2n}}
= &\ \left(\prod_{k=1}^{n-2}\prod_{r=0}^{k-1}
\balpha_{k,aq^{n-k-1+2r}} \balpha_{j,aq^{n+k-1-2r}}\right)
\left(\prod_{r=0}^{n-2}\balpha_{n-1,aq^{2r}}\balpha_{n,aq^{2+2r}}\right).
\end{align*}

\end{enumerit}
\end{lem}

\subl{} We can now prove Proposition \ref{gammaq}.

\begin{pf} Assume first that $\lie g$ is of type $A_n$. We claim
that the assignment $\chi_a\to\bomega_{1,a}$ defines a
homomorphism $\tau: \Xi_q\to \cal{P}_q/\cal{Q}_q$. For this, it is
enough to check that for all $a\in\bc^\times$,
\begin{equation}\label{trivan} \prod_{r=0}^n\bomega_{1,aq^{n+1-2r}}\in\cal{Q}^+_q.\end{equation}
 By using  Lemma
\ref{ellfund}(i) repeatedly we see that
$$\bomega_{n,q}=\left(\prod_{r=0}^{n-1}\bomega_{1,q^{n-2r}}\right)(\varpi)^{-1},$$
for some $\varpi\in \cal{Q}_q^+.$  Hence to prove \eqref{trivan}
it suffices to observe from Lemma \ref{ellfund}(i)  that
$$\bomega_{1,q^{-n-1}}\bomega_n=\prod_{j=1}^n\balpha_{j,
q^{-n+j-2}}.$$

 To see that this is an isomorphism of groups,
consider first the homomorphism $\cal{P}_q\to\Xi_q$ given by
mapping $$\bomega_{i,a}\mapsto \sum_{r=0}^{i-1}\chi_{aq^{2r-i+1}}.$$
We claim that $\cal{Q}_q$ is in the kernel of this map.  For this
it is enough to prove that $\balpha_{i,a}$ is in the kernel. We
prove this by induction on $i$. If $i=1$, then the result follows
since
$\bomega_{2,aq}=\bomega_{1,aq^2}\bomega_{1,a}\balpha_{1,a}^{-1}$.
The inductive step is now easily completed using Lemma
\ref{ellfund}(i) again.
 Thus we have a  homomorphism $\cal{P}_q/\cal{Q}_q\to \Xi_q$ which is clearly an
 inverse of $\tau$ and we are done.

Assume next that $\lie g$ is of type $B_n$. Define a group
homomorphism $\cal{P}_q\to\Xi_q$ by extending
$$\bomega_{i,a}\mapsto\chi_{aq^{2n-2i-1}}+\chi_{a,q^{-2n+2i+1}},\ i<n,\
\ \bomega_{n,a}\mapsto\chi_a.$$ Using Lemma \ref{ellfund}(ii) we see
by an induction starting at $n$ that $\balpha_{j,a}$ is the kernel
of this map and hence we get a homomorphism from
$\cal{P}_q/\cal{Q}_q\to \Xi_q$. To see that this map is an
isomorphism it suffices to show as in the case of $A_n$ that the
assignment $\chi_a\mapsto\bomega_{n,a}$ defines a homomorphism
$\Xi_q\to\cal{P}_q/\cal{Q}_q$. For this, it is enough to show that
$$\bomega_{n,a}\bomega_{n,aq^{4n-2}}\in\cal{Q}_q^+,$$ which is
just the second statement in Lemma \ref{ellfund}(ii).

If $\lie g$ is of type $C_n$, then we show by using Lemma
\ref{ellfund}(i) that the map $\chi_a\mapsto \bomega_{1,a}$ gives an
isomorphism between $\Xi_q$ and $\cal{P}_q/\cal{Q}_q$. We omit the
details.
\end{pf}

\begin{rem} If $q=1$, note that $\Xi_1$ is isomorphic to the
group of functions from $\bc^\times\to \Gamma$ with finite
support (cf \cite{CM}).
\end{rem}

\section{ The $\ell$--weights of finite--dimensional representations}

\subl{}

 Given a $\bu$-module $V$ and $\mu=\sum_i\mu_i\omega_i\in
P$, set
\begin{equation*} V_\mu=\{ v\in V: K_i.v =q_i^{\mu_i}v ,\ \
\forall \ i\in I\}.\end{equation*} We say that $V$ is a module of
type 1 if
\begin{equation*} V=\bigoplus_{\mu\in P}V_\mu.\end{equation*}
 Analogous definitions hold for representations of $\bu^{fin}$. Recall from
\cite{L} that for every $\lambda\in P^+$, there exists a unique
(up to isomorphism) irreducible finite--dimensional representation
of $\bu^{fin}$ which we denote by $V(\lambda)$. Let $\cal C_q$
be the abelian category consisting of type 1 finite--dimensional
representations of $\bu$. We set $$\wt(V)= \{\mu\in P: V_\mu\ne
0\}, $$ and given $v\in V_\mu$ we set $\wt(v)=\mu$.

\subl{} \begin{defn} Let $V$ be a $\bu$--module.  We say that
 $\bvpi\in\cal{A}^n$
is an $\ell$--weight of $V$ if there exists a non--zero element
$v\in V$ such that $$(P_{i, r}-(\varpi_i)_r)^N v=0, \ \ \ N\equiv
N(i,r,v)\in\bz^+,$$ for all $i\in I$ and $r\in\bz^+$ and we call
$v$ an $\ell$--weight vector in $V$ with $\ell$--weight $\bvpi$.
Denote the subspace consisting of all $\ell$--weight vectors with
$\ell$--weight $\bvpi$ by $V_{\bvpi}$.
\end{defn}
\begin{rem} We shall see later in the section, (see Proposition \ref{pm}), that the
generalized eigenspaces for the action of the $P_{i,-r}$, $i\in
I$, $r\in\bz_+$ are actually determined uniquely by those of the
$P_{i,r}$, $r\in\bz_+$.
\end{rem}
\subl{} The following lemma is trivially established.

\begin{lem}
Let $V\in\cal{C}_q$. We have
$$V=\opl_{\bvpi\in\cal{A}^n}V_{\bvpi},\ \ \ \
V_\mu=\opl_{\bvpi\in\cal{A}^n}V_{\bvpi} \cap V_\mu
.$$\hfill\qedsymbol\end{lem}
 Denote by $\wtl(V)$ the set of $\ell$--weights of $V$ and
define $\wtl(v)$ in the obvious way. It is obvious from the
definition of $\ell$--weights that any morphism between objects of
$\cal{C}_q$ preserves $\ell$--weight spaces.

\subl{}  We now study the behavior of $\ell$--weights
under tensor products. This is essentially the same proof given in
\cite{FR}, we include it here for completeness.
\begin{lem}\label{ltensor}  Let $V_r\in\cal{C}_q$, $r=1,2$ and  let $v_{j,r}$,
$1\le j\le \text{dim}(V_r)$ be a basis of $V_r$, $r=1,2$ such that
$\wt(v_{j,r})$ and $\wtl (v_{j,r})$ are defined and assume that if
$j<j'$ then $\wt(v_{j,r})-\wt(v_{j',r})\in Q^+$. Then the
$\ell$--weight vectors of $V_1\otimes V_2$ are of the form
$$v_{j,1}\otimes v_{j',2} +\left(\text{terms in}\ \oplus
(V_1)_{\nu_1}\otimes (V_2)_{\nu_2}\right),$$ where the direct sum
is over $\nu_r\in\wt(V_r)$, $r=1,2$ and $\nu_2-\wt(v_2)\in Q^+$.
The corresponding $\ell$--weight is $\wtl(v_{j,1})\wtl
(v_{j',2})$. In particular, $$\wtl(V_1\otimes
V_2)=\wtl(V_1)\wtl(V_2).$$
\end{lem}
\begin{pf} It is easy to
 see from Proposition \ref{comultip} that the matrices of the
action of  $P_{i, s}$, $i\in I$, $s\in\bz_+$ on $V_1\otimes V_2$
with respect to the basis $v_{j,1}\otimes v_{j',2}$ are
simultaneously upper triangular with diagonal entries given by
$\wtl(v_{j,1})\wtl (v_{j',2})$. The result is now immediate.
\end{pf}

\subl{}We need several results on irreducible representations of
$\bu$. We begin with the definition of an $\ell$--highest weight
module.
\begin{defn} We say that a $\bu$--module $V\in\cal{C}_q$ is
$\ell$--highest weight with $\ell$--highest weight
$\bvpi\in\cal{P}_q$  if there exists a non--zero vector $0\ne v\in
V$ such that $V=\bu v$ and, \begin{equation}\label{fdrel}
x_{i,r}^+v =0,\ \ P^\pm_i(u)v=(\bvpi)^\pm_iv,\ \ K_i^{\pm
1}v=q^{\pm\wt\bvpi(\check{\alpha}_i)} v,\ \
(x_{i,r}^-)^{\wt\bvpi(\check{\alpha}_i)+1} v=0,\end{equation} for
all $i\in I$, $r\in\bz$. The element $v$ is called the
$\ell$--highest weight vector.
\end{defn}
The following lemma is standard.
\begin{lem} Any $\ell$--highest weight module has a unique
irreducible quotient which is also a highest weight module with
the same highest weight.\end{lem}

 \subl{} The following was proved in \cite{CPqa},\cite{CPbanff} .
\begin{thm}\label{class}\hfill
\begin{enumerit}
\item[(i)] Any irreducible module in $\cal{C}_q$ is
$\ell$--highest weight.
\item[(ii)] There exists a bijective correspondence  between
elements of $\cal{P}_q^+$ and isomorphism classes of  irreducible
modules in $\cal{C}_q$.\hfill\qedsymbol
\end{enumerit}
\end{thm}
\begin{cor} Let $V\in\cal{C}_q$ be a highest weight module with
highest weight $\bvpi$. Then
$\bvpi\in\cal{P}_q^+$.\hfill\qedsymbol\end{cor}

\subl{} Given $\bomega\in\cal{P}_q^+$, let $V(\bomega)\in\cal{C}_q$
be an element in the corresponding isomorphism class, and let
$v_\bomega$ be the $\ell$--highest weight vector. We note the
following simple consequence of Theorem \ref{class}.
\begin{lem}
\label{classcon} Let $\bomega\in\cal{P}_q^+$.\hfill
\begin{enumerit}
\item[(i)]We have $V(\bomega)_{\wt\omega}=\bc v_\bomega.$
\item[(ii)] If
 $v\in V(\bomega)$ is such that
$x_{i,k}^+v=0$ for all $i\in I$, $k\in\bz$, then, $v=cv_\bomega$
for some $c\in\bc$.
\item[(iii)] Let $V\in\cal{C}_q$, $\bvpi=(\varpi_1,\cdots ,\varpi_n)\in\cal{C}_q$ and assume that $v\in V_\bvpi$, $v\ne 0$, is such that
$x_{i,k}^+v=0$ for some $i\in I$ and all $k\in\bz$. Then,
$\varpi_i\in\bc[u]$.\end{enumerit} \hfill\qedsymbol\end{lem}

\subl{} We now consider the case when $\lie g=sl_2$. Given $m>0$,
let $\bomega_a(m)\in\cal{P}_q^+$ be the polynomial $$\bomega_a(m)
= (1-aq^{m-1}u)(1-aq^{m-3}u)\cdots (1-aq^{-m+1}u),$$ and set
$\bomega_a(0)=1$. The following result was proved in \cite{CPqweyl}.

\begin{thm}\label{weylmsl2}
Suppose that  $V\in\cal C_q$ is an $\ell$-highest weight module
with highest weight $\bomega= \prod_{r=1}^k(1-a_ru)$ where
$a_r\in\bc^\times$ is such that $a_r/a_{r'}\ne q^2$ if $r'<r$.
Then $V$ is a quotient of $V(\bomega_{a_1}(1))\otimes\cdots
\otimes V(\bomega_{a_k}(1))$.\hfill\qedsymbol

\end{thm}

\subl{}
\begin{prop}\label{sl2class} \hfill
\begin{enumerit}
\item[(i)] There exists an isomorphism of $\bu^{fin}$--modules
$$V(\bomega_a(m))\cong V(m\omega_1).$$
\item[(ii)] The
eigenvalues of $P_1^\pm(u)$ on $V(\bomega_a(m))$ are
$\bvpi_{a,r}(m)^\pm$, $0\le r\le m$, where $$\bvpi_{a,r}(m)^\pm=
\bomega_{aq^{-r}}(m-r)^\pm(\bomega_{aq^{m-r+2}}(r)^\pm)^{-1}.$$ In
particular $\bvpi_{a,m}=\bvpi_{aq^2}(m)^{-1}$.
\item[(iii)] Any irreducible $\bu$--module is isomorphic to a tensor
product $V(\bomega_{a_1}(m_1))\otimes \cdots\otimes
V(\bomega_{a_r}(m_r))$ where $$a_k/a_s\ne q^{\pm(m_k+m_s-2p)},\ \
0\le p<\min\{m_k,m_s\},$$ for some $a_1,\cdots ,a_r\in\bc^\times$
and $m_1,\cdots ,m_r\in\bz^+$.
\end{enumerit}
\end{prop}

\begin{pf} Parts (i) and (iii) were proved in \cite{CPqa}. Part
(ii) was proved in a slightly different form in \cite{FR}. We
include a proof here for the reader's convenience.

We proceed by induction on $m$. If $m=1$, let  $v_0,v_1$ be  the
basis for $V(\bomega_a(1))$, where $v_0$ is the $\ell$--highest
weight vector. It is now a simple computation to check, using the
formulas in \cite{CPqa} to check that the eigenvalues of
$P_1(u)^\pm$ on $v_1$ are $\bvpi_{aq^2}(1)^{-1}=(1-aq^2u)^{-1}$.
 Assume now that we know the
result for all $s<m$. By \cite{CPqa}, we know that there exists a
short exact sequence, $$0\to V(\bomega_{aq}(m-2))\to
V(\bomega_{a}(m-1))\otimes V(\bomega_{aq^{-m}}(1)) \to
V(\bomega_{aq^{-1}}(m))\to 0.$$ Let $v_0,v_1$ be a basis for
$V(\bomega_{aq^m}(1))$ and $w_0,\cdots, w_{m-1}$ a basis for
$V(\bomega_{a}(m-1))$. Then, the elements $w_j\otimes v_0$, $0\le
j\le m-1$ and the element  $w_{m-1}\otimes v_1$ all must have
non--zero projection onto $V(\bomega_{aq}(m))$. For otherwise,
applying $x_0^+$ repeatedly we find that $w_0\otimes v_0\in
V(\bomega_{aq^{-1}}(m-2))$ which is impossible. On the other hand,
using the formulas in Proposition \ref{comultip} we see that
$$P(u) (w_j\otimes v_0)=
(\varpi^j_{a,m-1}\varpi^0_{aq^{-m},1})(w_j\otimes v_0)=
(\varpi^j_{aq^{-1},m})(w_j\otimes v_0),\ \ 0\le j\le j-1$$ and
$$P(u)(w_{m-1}\otimes v_1)=\frac{1}{\bomega_{aq}}(w_{m-1}\otimes
v_1) =(\varpi^m_{aq^{-1},m})(w_{m-1}\otimes v_0) .$$ The result
follows.
\end{pf}

\subl{} We now assume  that $\lie g$ is an arbitrary simple Lie
algebra.
\begin{prop}\label{pm} Let $V\in\cal{C}_q$.
\begin{enumerit}
\item[(i)] Suppose that  $0\ne v\in V_\bvpi$ is such that
$x_{i,k}^+v=0$ for some $i\in I$ and all  $k\in\bz$. Then
$(\bvpi)_i\in\bc[u]$.
\item[(ii)] For all $\bvpi\in\wt_\ell(V)$, we have $\bvpi\in\cal{P}_q$ and  $V_\bvpi\subset
V_{\wt\bvpi}$. In particular, $$\wt(V)=\{\wt(\bvpi):
\bvpi\in\wt_\ell(V)\}.$$
\item[(iii)] The eigenvalues of the elements $P_{i,-r}$, $i\in I$,
$r\in\bz$, on $V_\bvpi$ are given by $\bvpi^-\in\cal{P}_q$.
\end{enumerit}
\end{prop}

\begin{pf} It clearly suffices to prove the proposition when $\lie
g=sl_2$. Since $V$ is finite--dimensional it has a Jordan--Holder
series and hence we may assume without loss of generality that $V$
is irreducible, say $V=V(\bomega)$ for some
$\omega\in\cal{P}_q^+$. The proof of (i) is now immediate since
Corollary \ref{class}(i) implies that $v=cv_\bomega$ for some
$c\in\bc^\times$. To prove the other parts notice that by
Proposition \ref{sl2class}(ii) and Lemma\ref{ltensor}, it suffices
to consider the case of the representations $V(\bomega_a(r))$,
$r\in\bz_+$, $a\in\bc^\times$. But this is exactly part (iii) of
Proposition \ref{sl2class}.
\end{pf}

\subl{} Let $\bz[\cal{P}_q]$ be the integral group ring over
$\cal{P}_q$ and let $e(\bvpi), \bvpi\in\cal P_q$ be a basis of the
group ring .

\begin{defn}
 Given $V\in\cal{C}_q$, let
$ch_\ell(V)\in\bz[\cal{P}_q]$ be defined by
$$ch_\ell(V)=\sum_{\bvpi\in\cal{P}_q}\dim(V_\bvpi)\ e(\bvpi).$$
\end{defn}
 In \cite{FR} it was proved that the
$q$--character of $V$, which was defined using the $R$--matrix, is
just $ch_\ell(V)$. It is quite clear, as observed in \cite{FR},
that $ch_\ell$ is additive and by Lemma \ref{ltensor}
multiplicative.

\section{Braid group invariance of $\ell$--weights of fundamental
representations} Throughout this section we assume that $\lie g$
is of classical type.  The representations $V(\bomega_{i,a})$ ,
$i\in I$, $a\in\bc^\times$ are called the fundamental
$\ell$--highest weight representations. The main result of this
section is the following theorem. Recall from Section 1, the
subsets $W_\lambda\subset W$ defined for elements $\lambda\in
P^+$.

\subl{}
\begin{thm}\label{wconjfund}  Let $i\in I$, $a\in\bc^\times$ and assume that
$\bvpi\in\wt_\ell(V(\bomega_{i,a}))$ is such that
$\wt(\bvpi)=\lambda\in P^+$.
\begin{enumerit}
\item[(i)] Let $w'=s_jw\in W_\lambda$ for some  $j\in I$ with $\ell(s_jw)=\ell(w)+1$. Then,
$(T_{w}\bvpi)_j\in\bc[u]$.  In particular,
$$\bvpi\in\bomega_{i,a}\cal{Q}_q^-\implies T_{w'}
\bvpi\in\bomega_{i,a}\cal{Q}_q^-.$$
\item[(ii)] For all $w\in W_\lambda$ we have
$$\dim(V(\bomega_{i,a})_\bvpi)=\dim(V(\bomega_{i,a})_{T_w\bvpi}),$$
and
$$T_w(\wt_\ell(V(\bomega_{i,a})_\lambda))=\wt_\ell(V(\bomega_{i,a})_{w\lambda}).$$
\item[(iii)] Suppose that $\bvpi\ne \bomega_{i,a}$.  There exists
$\bvpi'=(\varpi_1',\cdots
,\varpi_n')\in\wt_\ell(V(\bomega_{i,a}))$ and $j\in I$ with
$\varpi_j'=(1-cu)(1-c'u)$, and
\begin{equation}\label{indom} \bvpi=\bvpi'(\balpha_{j,c})^{-1}.\end{equation}
If $c'\ne cq^{-2}$ then $\bvpi'\ (\balpha_{j,c'})^{-1}\in\wtl
V(\bomega_{i,a})$  and  if $c=c'$ then
$\dim(V(\bomega_{i,a}))_{\bvpi}\ge 2$.

\end{enumerit}
\end{thm}

\subl{}
\begin{cor} \label{cone} We have
$\wtl(V(\bomega_{i,a}))\subset\bomega_{i,a}\cal{Q}_q^-.$
\end{cor}
\begin{pf} By part (i) of the theorem, it suffices to prove the
corollary  for $\bvpi\in\wtl(V(\bomega_{i,a}))$ with
$\wt(\bvpi)=\lambda\in P^+$.  We proceed by induction on
$\rm{ht}(\omega_i-\lambda)$, with induction obviously beginning
when $\omega_i=\lambda$. Let $\bvpi'$ be as in part (iii) of the
theorem, so that  $\wt(\bvpi')=\lambda+\alpha_j$. Choose $w\in W$
and $\mu\in P^+$ with $w\mu=\lambda+\alpha_j$. Then
$\mu\ge\lambda$ and by part (ii), we have $$\bvpi'=T_w\bvpi''$$
for some $\bvpi''\in\wtl(V(\bomega_{i,a}))$ with
$\wt(\bvpi'')=\mu$. By the induction hypothesis
$\bvpi''\in\bomega_{i,a}\cal{Q}_q^-$ and hence by part (i) again
$\bvpi'\in\bomega_{i,a}\cal{Q}_q^-.$ The corollary is proved.
\end{pf}

\begin{rem} The corollary was proved \cite{FM} by combinatorial methods for
all simple Lie algebras. On the other hand, Theorem
\ref{wconjfund} is not true for the exceptional algebras since
there is no suitable analog of Lemma 1.3 available for those
algebras. In particular, it follows that Theorem \ref{wconjfund}
is stronger than the Corollary. The theorem is also obviously
false for an arbitrary irreducible representation since the
evaluation representations of $\bu(\lie sl_2)$ of Section 4 are
counterexamples.
\end{rem}

\subl{} \begin{cor} We have  $$
ch_\ell(V(\bomega_{i,a}))=\sum_{\bvpi\in \cal{P}_q,\wt\bvpi=\mu\in
P^+} \dim (V_{\bvpi})\left(\sum_{w\in W_\mu} \ e(T_w\bvpi)\right).$$
\end{cor}
We remark here that the preceding results are analogous to the
following well--known result for finite--dimensional
representations $V$ of simple Lie algebras: for all $\mu\in P$
with $V_\mu\ne 0$, and $w\in W$, we have
$$\dim (V_\mu)=\dim (V_{w\mu}).$$

\subl{} Before proving Theorem \ref{wconjfund} we note some
consequences of it in computing $q$--characters (or
$\ell$--characters) of the fundamental representations, see also
\cite{Nak}. Thus, suppose that $i\in I$ is such that the
$\bu^{fin}$--representation $V(\omega_i)$ is miniscule.  In that
case there exist no weights $\mu\in P^+$ such that $\mu\le
\omega_i$ and it is known that $V(\bomega_{i,a})\cong V(\omega_i)$
as $\bu^{fin}$--modules. Hence, the $q$--character of the
fundamental representation is of the
form,$$ch_\ell(V(\bomega_{i,a}))=\sum_{w\in W_{\omega_i}}
e(T_w\bomega_{i,a}).$$

\subl{} Suppose next that $\lie g$ is of type $D_n$, and that $V=V(\bomega_{2,a})$. Let $\bvpi_j$ be defined by
\begin{alignat*}{5}
&\bvpi_j  && = (\bomega_{j-1,aq^{j+1}})^{-1} \bomega_{j-1,aq^{2n-j-3}}\bomega_{j,aq^{j}}(\bomega_{j,aq^{2n-j-2}})^{-1},&& \text{ for } 1\leq j \leq n-2,\\
&&& = \bomega_{j,aq^{n-3}}(\bomega_{j,aq^{n+1}})^{-1}, && \text{ for } j=n-1,n.
\end{alignat*}
We claim that
\begin{equation}\label{adjoint}
ch_\ell(V)=\sum_{w\in W_{\omega_2}} e(T_w\bomega_{2,a})\
+ \ \sum_{j\neq n-2} e(\bvpi_j) +2 e(\bvpi_{n-2}).
\end{equation}
 To prove the claim, set $$w_j =
s_{j-1}\cdots s_1s_{j+1}\cdots s_{n-2} s_{n}s_{n-1} s_{n-2}\cdots
s_1,$$ and observe that  $w_j\in W_{\omega_2}$ with  $w_j\omega_2
= \alpha_j$. A straightforward calculation shows that,
\begin{alignat*}{5}
&T_{w_j}\bomega_{2,a} && = (\bomega_{j-1,aq^{j+1}})^{-1} \bomega_{j,aq^{j}} \bomega_{j,aq^{2n-4-j}} (\bomega_{j+1,aq^{2n-3-j}})^{-1},&& \ \text{ if } j\leq n-3,\\
&&&=(\bomega_{n-3,aq^{n-1}})^{-1} (\bomega_{n-2,aq^{n-2}})^{2}(\bomega_{n-1,aq^{n-1}}\bomega_{n,aq^{n-1}})^{-1},&&  \ \text{ if } j= n-2,\\
&&& = (\bomega_{n-2,aq^{n}})^{-1} \bomega_{j,aq^{n-1}}\bomega_{j,aq^{n-3}},&&\ \text{ if } j=n-1,n.
\end{alignat*}
It follows from Theorem \ref{wconjfund} that $$\wtl(V) =
\{T_w\bomega_{2,a}:w\in W_{\omega_2}\}\cup\{\bvpi_j:1\le j\le n\},$$ and that
$$\dim(V_{\bvpi_{j}}) \ge 1,\ \ j\ne n-2,\ \
\dim(V_{\bvpi_{n-2}})>1.$$ On the other hand, we  know from
 \cite{CPbanff} that $V\cong V(\omega_2)\oplus\bc$ as
$\bu^{fin}$--modules. In particular, $\dim(V_0)=n+1$. Hence it
follows that $$\dim(V_{\bvpi_{j}}) = 1,\ \ j\ne n-2,\ \
\dim(V_{\bvpi_{n-2}})=2$$ and our claim is proved.

\subl{} We shall need the following result which holds for all
simple Lie algebras $\lie g$.
\begin{prop}\label{extr1} Let   $V\in\cal{C}_q$, $\bvpi=(\varpi_1,\cdots ,\varpi_n)\in\cal{P}_q$ and  suppose that there exists  $v\in
V_\bvpi$, $v\ne 0$, and $j\in I$ such that $$x_{j,s}^+ v =0,\ \ \
\forall \ s\in\bz.$$
\begin{enumerit}
\item[(i)] Then,
$\varpi_j\in\bc[u]$ is of degree $\wt(\bvpi)(\check\alpha_j)$ and
$$(x_{j,0}^-)^{\wt(\bvpi)(\check\alpha_j)}v\in V_{T_j\bvpi}\ \ .$$
\item[(ii)] Write $\varpi_j$ as a product of polynomials
$\prod_{r=1}^{k}\bomega_{a_r}(m_r)$ as in Proposition \ref{sl2class}.
Then $\bvpi(\balpha_{j,a_rq^{m_r-1}})^{-1}\in\wtl(V)$ for all $1\le r\le k$. Furthermore,
 for all $s\in\bz$,  we have
$$x_{j,s}^-v\in \sum_{r=1}^k \sum_{l=0}^{m_r-1}V_{\bvpi\balpha_{j,a_rq^{m_r-1-2l}}^{-1}},$$
and
 $$\dim V_{\bvpi\balpha_{j,a_rq^{m_r-1}}^{-1}}\ge \#\{1\le s\le k:a_r=a_s\}.$$

\end{enumerit}
 Analogous
statements hold if $x_{j,s}^-v=0$ for all  $s\in \bz$.
\end{prop}

\begin{pf} By Lemma \ref{classcon}  and Proposition \ref{pm}, we see that $\varpi_j\in\bc[u]$ and
$\wt(\bvpi)(\check\alpha_j)={\rm deg} \varpi_j$. Further,
$(x_{j,0}^-)^{\wt(\bvpi)(\check\alpha_j)}v\ne 0$. Suppose first
that $v$ is actually an eigenvector  with eigenvalue $\bvpi$. It
follows from Lemma \ref{classcon}(iii) and   Proposition
\ref{sl2class}(iii), that
$(x_{j,0}^-)^{\wt(\bvpi)(\check\alpha_j)}v$ is an eigenvector for
$P_j^+(u)$ with eigenvalue $(\varpi_j(q_j^2u))^{-1}$. To compute
the eigenvalues corresponding to $P_i^+(u)$ it suffices in view of
\eqref{hp}  to compute the eigenvalues for $h_{i,r}$ for all $i\in
I$ and $r\in\bz$, $r\ne 0$.  To simplify our notation, we set for
all $i\in I$, $r\in\bz$, $r>0$,  $$\tilde h_{i,r} =
-\frac{q^r}{[r]_{q_i}}h_{i,r},\ \ k=\wt(\bvpi)(\check\alpha_j).$$
Using the relations in $\bu$, we see that
\begin{equation}\label{comm}[\tilde h_{i,r},(x_{j,s}^{-})^k] =
\frac{[ra_{ij}]_{j}[r]_{j}}{[2r]_{j}[r]_{i}}[\tilde
h_{j,r},(x^-_{j,s})^k] =
\frac{[ra_{ij}]_{j}}{(q_j^r+q_j^{-r})[r]_{i}}[\tilde
h_{j,r},(x^-_{j,s})^k],\end{equation} which gives,
\begin{align*}
\tilde h_{i,r}(x_{j,s}^-)^k = (x_{j,s}^-)^k\tilde h_{i,r} +
[\tilde h_{i,r},(x_{j,s}^{-})^k] = (x_{j,s}^-)^k \tilde h_{i,r} +
\frac{[ra_{ij}]_{j}}{(q_j^r+q_j^{-r})[r]_{i}}[\tilde
h_{j,r},(x^-_{j,s})^k].
\end{align*}
 Writing
$\ln(\varpi_i)=\sum_{r\ge 1}\varpi_{i,r}u^r$, for some
$\varpi_{i,r}\in\bc$, we find from \eqref{hp} that,
 $$\tilde h_{i,r}v =
\varpi_{i,r}v,\ \ i\in I, \ r\in\bz, r>0.$$ On the other hand, we
have already observed that  $$ \tilde h_{j,r}(x_{j,s}^-)^kv= (
\ln(\varpi_j(q_j^2u)))_r(x_{j,s}^-)^kv.$$ Since  $\varpi_j$ is a
polynomial of degree $k$, write $$\varpi_j(u) = \prod_{t=1}^k
(1-b_tu)$$ for some $b_1,\cdots ,b_k\in\bc^\times$. This means
that  $$(\ln(\varpi_j(q_j^2u)))_r = (\sum_{t=1}^k
b_t^rq_j^{2r})/r.$$ It is now a simple checking to see that
\begin{align*}
\tilde h_{i,r}(x_{j,s}^-)^kv = \left(\varpi_{i,r} +
\frac{q_j^r[ra_{ij}]_{j}}{r[r]_{j}}\sum_{t=1}^k b_t^r
\right)(x_{j,s}^-)^kv
\end{align*}
for all $i\in I$ and $r\in\bz$, $r>0$. Using \eqref{hp} again we
see that, $$(x_{j,0}^-)^{\wt(\bvpi)(\check\alpha_j)}v\in
V_{T_j\bvpi}.$$ For the second statement, observe first that, by
Theorem \ref{weylmsl2} and Proposition \ref{sl2class}, $x_{j,s}^-v$ is a sum of eigenvectors
for $P_j^+(u)$ with eigenvalues
$(1-a_rq^{m_r-1-2l+2}u)^{-1}(1-a_rq^{m_r-1-2l}u)^{-1}\varpi_j$ where
$\varpi_j=\prod_{r=1}^k\bomega_{a_r}(m_r)$ and $l=0,1,\cdots, m_r-1$. But now, a calculation
identical to the preceding one gives the result. The statement on
the dimensions follows from Lemma \ref{ltensor}, Theorem \ref{weylmsl2} and Proposition \ref{sl2class}.

It remains to consider the case when $v$ is a generalized
eigenvector for the action of the $P_{i}(u)$. Clearly, we can
choose a Jordan basis $v_1,\cdots ,v_m$ of $\bu(0)v\subset
V_\bvpi$ simultaneously for the action of the $\tilde h_{i,k}$,
i.e., $$\tilde h_{i,r}v_t\in \oplus_{t'< t}\bc v_{t'}
+\varpi_{i,r}v_t$$ for all $1\leq t\leq m$. We proceed by
induction on $t$, the case $t=1$ is dealt with above. Let $V_t$ be
the $\bu^j$--module generated by $v_t$. Then $V_{t}\subset
V_{t+1}$ and the image of $v_{t+1}$ in the quotient
$V_{t+1}/V_{t}$ is an $\ell$--highest weight vector for $\bu^j$.
In particular, it follows  that $$\tilde
h_{j,r}(x_{j,s}^-)^kv_t\in \oplus_{t'<t}\bc
(x_{j,s}^-)^kv_{t'}+\tilde\varpi_{j,r}(x_{j,s}^-)^kv_t$$ where
$\tilde\varpi_{j,r}=(\ln(\varpi_j(q_j^2u)))_r$. Using the
inductive hypothesis and \eqref{comm} we get
\begin{align*}
\tilde h_{i,r}(x_{j,s}^-)^kv_t \in \oplus_{t'<t}\bc (x_{j,s}^-)^kv_{t'} + \left(\varpi_{i,r} +
\frac{q_j^r[ra_{ij}]_{j}}{r[r]_{j}}\sum_{t=1}^k b_t^r \right)(x_{j,s}^-)^kv_t
\end{align*}
It follows that $(x_{j,0}^-)^kv_t\in V_{T_j\bvpi}$. The second
statement is proved similarly.
\end{pf}

\subl{} We now prove Theorem \ref{wconjfund}.

\begin{pf} Notice that by Lemma \ref{domfund} we can assume that
either $\lambda=0$ or $\lambda=\omega_r$ for some $r\in I$. Assume
first that $\lambda=\omega_r$ and let $w'=s_jw$ where
$\ell(w')=\ell(w)+1$. It follows that  $w^{-1}\alpha_j\in Q^+$ and
hence $(w\omega_r,\alpha_j)=(\omega_r,(w)^{-1}\alpha_j) > 0$. In
particular, this means that $(w)^{-1}\alpha_j-\alpha_r\in Q^+$ and
hence by Lemma \ref{domfund} we see that
$$\omega_i-(\omega_r+(w)^{-1}\alpha_j)\notin Q^+.$$ This in turn
implies that
$$\omega_r+(w)^{-1}\alpha_j\notin\wt(V(\bomega_{i,a})),$$ or
equivalently that $w\omega_r+\alpha_j\notin\wt(V(\bomega_{i,a}))$.
Thus,
\begin{equation}\label{jhw}
 x_{j,s}^+V(\bomega_{i,a})_{T_{w}\bomega}=0,\ \ \forall\
s\in\bz.
\end{equation}
A similar argument proves that
$$x_{j,0}^-V(\bomega_{i,a})_{T_{w'}\bomega}= 0.$$
 The first statement in part (i) is immediate from Lemma
\ref{classcon} while the second follows from equation
\eqref{braid2} and the fact that $T_{w'}=T_jT_{w}$. To prove (ii),
note that   Proposition \ref{extr1} implies that the map
$(x_{j,0}^-)^{(w\omega_r,\alpha_j)}:
V(\bomega_{i,a})_{T_{w}\bomega}\to
V(\bomega_{i,a})_{T_{w'}\bomega}$ is an isomorphism of vector
spaces. Since, $\dim V(\bomega_{i,a})_{\lambda}=\dim
V(\bomega_{i,a})_{w\lambda}$ th second statement in (ii) follows
as well. If $\lambda=0$ there is nothing to prove in (i) and (ii)
since $W_\lambda=\{e\}$.

 To prove (iii), choose a non--zero element
$v\in V(\bomega_{i,a})_\bvpi$,
 $j\in I$, and $s\in\bz$ such that $v=x_{j,s}^-v'$. Then $v'\in
V(\bomega_{i,a})_{\lambda+\alpha_j}\ne 0$ and by Lemma
\ref{domfund} we see that
$$x_{j,s}^+V(\bomega_{i,a})_{\lambda+\alpha_j}=0.$$
 Since
$(\lambda+\alpha_j)(\check\alpha_j)=2$,  it follows from Lemma
\ref{classcon} that $\varpi_j$ is a polynomial of degree two for
all $\bvpi'=(\varpi_1,\cdots, \varpi_n)\in\wtl(V(\bomega_{i,a}))$
with $\wt(\bvpi')=\lambda+\alpha_j$. Furthermore, since $$v'\in
\oplus_{\wt\bvpi'=\lambda+\alpha_j} V(\bomega_{i,a})_{\bvpi'},$$
it follows  from  Proposition \ref{extr1} that
$$\bvpi=\bvpi'(\balpha_{j,c})^{-1},$$ for some
$\bvpi'\in\cal{P}_q$ with $\wt\bvpi'=\lambda+\alpha_j$ and
$c\in\bc^\times$ satisfying $\varpi'_j(c^{-1})=0$.

To prove the final statement in (iii), let $v''\in
V(\bomega_{i,a})_{\bvpi'}$, where $\wt\bvpi'=\lambda+\alpha_j$ and
$\varpi_j'=(1-cu)(1-c'u)$. Since $x_{j,s}^+v''=0$ it follows from
Proposition \ref{sl2class} and Proposition \ref{extr1} that if
$c\ne c'q^2$, $V(\bomega_{i,a})_{\bvpi'(\balpha_{j,c})^{-1}}\ne 0$
and $ V(\bomega_{i,a})_{\bvpi'(\balpha_{j,c'})^{-1}}\ne 0$. The
other statement is proved similarly.
\end{pf}

\section{On the tensor product structure of Weyl modules}
\subl{} In this section we establish the conjecture in
\cite{CPqweyl} on the structure of finite--dimensional Weyl
modules using some deep results of Nakajima. This allows us to
prove the following generalization of Corollary \ref{cone} and
\cite[Theorem 4.1]{FM}.
\begin{thm}\label{gencone} Let $V\in\cal{C}_q$ be an $\ell$--highest weight
representation with highest weight $\bomega\in\cal{P}_q^+$. We
have $$\wtl(V)\subset\bomega\cal{Q}_q^-.$$\end{thm}

\subl{} We begin by recalling the definition and some results on
Weyl modules. Thus, let $W(\bomega)$ be the $\bu$--module
generated by an element $v_\bomega$ satisfying the relations:
\begin{equation}\label{fdrel2}
x_{i,r}^+v_\bomega =0,\ \ P^\pm_i(u)v_\bomega=(\bomega)^\pm_iv_\bomega,\ \
K_i^{\pm 1}v_\bomega=q^{\pm\wt(\bomega)(\check\alpha_i)} v_\bomega,\ \
(x_{i,r}^-)^{\wt(\bomega)(\check\alpha_i)+1} v_\bomega=0,\end{equation} The
following result was proved in \cite{CPweyl}.
\begin{prop}\label{weylfindim}\hfill
\begin{enumerit}
\item[(i)] For all $\bomega\in\cal{P}_q^+$, we have
$W(\bomega)\in\cal{C}_q$.
\item[(ii)] Any $\ell$--highest weight representation in
$\cal{C}_q$ is a quotient of $W(\bomega)$ for some
$\bomega\in\cal{P}_q^+$.

\end{enumerit}\hfill\qedsymbol
\end{prop}

\subl{} We shall also need, the following result proved in
\cite{AK},\cite{C1},\cite{VV}.
\begin{thm}\label{cyc} Let $k\in\bz$, $k\ge 1$ and let $i_1,\cdots ,i_k\in I$, $a_1,\cdots ,a_k\in\bc^\times$. The $\bu$--module
$V(\bomega_{i_1,a_1})\otimes\cdots\otimes V(\bomega_{i_k,a_k})$ is
cyclic on the tensor product of the $\ell$--highest weight vectors
if for all $1\le s'<s\le k$ we have $a_s\ne a_{s'} q^r$ for any
$r\in\bz$, $r>0$.\hfill\qedsymbol\end{thm} The following corollary
is now  immediate.
\begin{cor}\label{fundquot}  There exist $i_1,\cdots ,i_k\in I$, $a_1,\cdots
,a_k\in\bc^\times$ such that
$V(\bomega_{i_1,a_1})\otimes\cdots\otimes V(\bomega_{i_k,a_k})$ is
a quotient of $W(\bomega)$ where,
$$\bomega_j=\prod_{k:i_k=j}(1-a_ku).$$ In particular, $$\dim
W(\bomega)\ge \prod_k
\dim(V(\bomega_{i_k,a_k})).$$\hfill\qedsymbol\end{cor}

\subl{} As a result, to prove Theorem \ref{gencone}, it now
suffices to prove that the $\ell$--weights of $W(\bomega)$ lie in
$\bomega\cal{Q}_q^-$. This is immediate from Lemma \ref{ltensor},
\cite[Theorem 4.1]{FM} and the following Theorem.

\begin{thm}\label{cpconj}  The module $W(\bomega)$ is isomorphic to a tensor
product of fundamental representations. In particular if
$i_1,\cdots ,i_k\in I$, $a_1,\cdots ,a_k\in\bc^\times$ are  such
that $\bv=V(\bomega_{i_1,a_1})\otimes\cdots\otimes
V(\bomega_{i_k,a_k})$ is cyclic on the tensor product of highest
weight vectors, then $\bv\cong W(\bomega)$ where $\bomega$ defined
as in Corollary \ref{fundquot}.
\end{thm}

\begin{rem} This theorem was conjectured in \cite{CPqweyl} where
it was proved when $\lie g=sl_2$. In the general case, the proof
we give  is deduced  easily from some deep results of Nakajima
\cite{Nak1} in the simply--laced case and of Beck and Nakajima in
 \cite{Nak2} in the general case.
\end{rem}

\begin{pf} Given $\lambda\in P^+$, let $\tilde{W}(\lambda)$ be the
$\bu$--module generated by an element $v_\lambda$ satisfying,
\begin{equation}\label{fdrel3}
x_{i,r}^+v_\lambda=0,\ \ K_i^{\pm 1}v_\lambda=q^{\pm\lambda_i}
v_\lambda,\ \ (x_{i,r}^-)^{\lambda_i+1} v_\lambda=0.\end{equation}
It is not hard to see that  $\tilde{W}(\lambda)$ can also be
regarded as a $\bu(0)$--module by right multiplication. Further,
it was shown in \cite{CPqweyl} that $$P_{i,\pm r}v_\lambda =0,\ \
i\in I, |r|>\lambda_i,\ \
P_{i,\lambda_i}P_{i,-\lambda_i}v_\lambda=v_\lambda.$$ The quotient
of $\bu(0)$ by the ideal $I(\lambda)$ generated by the elements,
$P_{i,r}$, $i\in I$, $|r|>\lambda_i$ and
$P_{i,\lambda_i}P_{i,-\lambda_i}-1$ can be identified with the
ring $\Sigma$ of symmetric functions in the variables $t_{i,r}$,
$i\in I$, $1\le r\le\lambda_i$ and $t_{i,\lambda_i}^{-1}$
\cite{CPqweyl},\cite{Nak1}. It was proved in \cite[Proposition
14.1.2] {Nak1} for the simply--laced case, and in \cite[Section 4]
{Nak2} for the general case, that $\bu(0)v_\lambda$ is free as a
module for $\Sigma$ of rank $$\prod_{i\in I}\dim
V(\bomega_{i,1})^{\otimes \lambda_i}.$$

Given $\bomega\in\cal{P}^+_q$, with $\wt\bomega=\lambda$, let
$\frak m_\bomega$ be the maximal  ideal in $\bu(0)$ generated by
the elements $P_{i,\pm r}-((\bomega^\pm)_i)_r$. The module
$W(\bomega)$ is clearly a quotient of $\tilde{W}(\lambda)$, and in
fact as vector spaces, we have $$W(\bomega)\cong
\tilde{W}(\lambda)/\frak m_\omega,$$ and hence $$\dim(W(\bomega))
=\prod_{i\in I}\dim V(\bomega_{i,1})^{\otimes \lambda_i}.$$The
result now follows from Corollary \ref{fundquot}. The second
statement is immediate.
\end{pf}

\section{Block Decomposition of $\cal{C}_q$}

It is by now well--known that the category $\cal{C}_q$ is not
semisimple. In this section we  use the results of the previous
section  to describe the   blocks in $\cal{C}_q$. To do this, we
redefine the notion of elliptic characters introduced in \cite{EM}
as elements of $\Xi_q$. This allows us to state and prove the main
result of \cite{EM} for generic $q$.

\subsection{} We begin by recalling the definition of the blocks of an abelian category.

Recall that that two objects $V_r\in\cal{C}_q$, $r=1,2$ are linked
if there does not exist a splitting of $\cal{C}_q$ into a direct
sum of abelian categories $\cal{C}_r$, such that $V_r\in\cal{C}_r$
for $r=1,2$. It is not hard to see that linking defines an
equivalence relation on $\cal{C}_q$ and a block is an equivalence
class of this relation.

\subsection{} We now define the notion of an elliptic character of
a representation.
\begin{defn}
The elliptic character of an irreducible representation
$V(\bomega)$ of $U_q$ is the element $\chi_{\bomega} =
\overline{\bomega}\in P_q/Q_q \cong \Xi_q$. A finite dimensional
representation $V$ of $\bu$ is said to have elliptic character
$\chi \in\Xi_q$ if every  irreducible constituent of $V$ has
elliptic character $\chi$.  Let   $\cal C_{\chi}$ be the
subcategory  of $\cal C_q$ consisting of representations with
elliptic character $\chi$.
\end{defn}

\begin{rem}
Recall in section \ref{gammaq} that we use additive notation when working with $\Xi_q$, rather than the multiplicative notation induced from $\cal P_q$.
\end{rem}

{\samepage
\begin{thm}\label{blockdec}\hfill
\begin{enumerit}
\item[(i)] Every indecomposable object in $\cal C_q$ has a well
defined elliptic character.
\item[(ii)] Any two simple objects in $\cal C_{\chi}$ are linked.
\item[(iii)] The categories $\cal C_{\chi}$ are the blocks of $\cal C_q$.
\end{enumerit}
\end{thm}}

We prove the theorem in the rest of the section.

\subsection{} The following   proposition plays an important  role in the proof of Theorem \ref{blockdec}.

\begin{prop}\label{tensor}\hfill
\begin{enumerit}
\item[(i)] For all $\bomega\in\cal{P}_q^+$, we have $W(\bomega)\in
\cal{C}_{\chi_\bomega}$. \item[(ii)] $\cal C_{\chi_1}\otimes \cal
C_{\chi_2} \subset \cal C_{\chi_1+\chi_2}$.
\end{enumerit}
\end{prop}
\begin{pf} To prove (i), note that if $V(\bomega')$ is a
constituent of $W(\bomega)$, then by Theorem \ref{gencone} we must
have $$\bomega'\in\bomega\cal{Q}_q^-.$$ It follows that $\chi_{\bomega}=\chi_{\bomega'}.$
The proof of the second part is similar. It is enough to prove that if $V(\omega_r)\in\cal
C_{\chi_r}$, $r=1,2$, then $V=V(\bomega_1)\otimes V(\bomega_2)\in
\cal C_{\chi_1+\chi_2}$. Suppose that $V(\bomega)$ is an
irreducible constituent of $V$. In particular $\bomega \in
\wt_\ell(V)\subset \wtl( V_1)\wtl( V_2)$. By Theorem \ref{cone} we
know that $\wtl( V_r)\subset\bomega_r\cal{Q}_q^-$ for $r=1,2$.
Hence $\bomega\in\bomega_1\bomega_2\cal{Q}_q^-$. Together with Proposition \ref{gammaq}, it
immediately implies that
$$\chi_{\bomega} = \overline{\bomega} = \overline{\bomega_1}
\overline{\bomega_2} = \chi_{\bomega_1}+\chi_{\bomega_2}.$$
\end{pf}

\subl{} We can now prove Theorem \ref{blockdec}(i) by methods
similar to those used in \cite{CM} for affine algebras.  Namely,
one proves the following:

\begin{lem}\label{nullext}\hfill
\begin{enumerit}
\item[(i)]
Let $U\in\cal C_\chi$. Let   $\bomega_0 \in \cal{P}_q^+$ be such
that $\chi\ne \chi_{\bomega_0}$. Then
${\rm{Ext}}_{\cal C_q}^1(U,V(\bomega_0)) = 0$.
\item[(ii)]
Assume that  $V_j\in\cal C_{\chi_j}$, $j=1,2$ and that $\chi_1\ne
\chi_2$. Then ${\rm{Ext}}_{\cal C_q}^1(V_1,V_2) = 0$.
\end{enumerit}
\end{lem}

\begin{pf} Since $\rm{Ext}^1_{\cal C_q}$ is an additive functor, to prove (i) it
suffices to consider the case when $U$ is indecomposable. Consider
an extension,
\begin{equation*}
0\to V(\bomega_0)\to V\to U \to 0
\end{equation*}
We prove  by induction on the length
of $U$ that the extension is trivial. Suppose first that
$U=V(\bomega)$ for some $\bomega\in\cal P_q^+$. Then $V_\bomega\ne 0$
and
 one of the following must hold,\begin{enumerit}
\item[(i)] $\wt\bomega<\wt\bomega_0$, or
\item[(ii)] $\wt{\bomega_0}-\wt\bomega \notin Q^+\backslash\{0\}$.
\end{enumerit}
We can always  assume (by taking duals if necessary) that we are
in case (ii). Since $\wt(V(\bomega_0))\subset\wt\bomega_0 - Q^+$,
it  follows that $$x_{i,k}^+ V_{\wt\bomega}=0,\ \ \forall\ \ i\in
I, k\in\bz. $$  Thus there exists an element $0\ne v\in
V_{\bomega}$ which is a common eigenvector for the action of
$P_i(u)$ with eigenvalue $\bomega$ and hence $\bu v$ is a quotient
of $W(\bomega)$. In particular $\bu v\in\cal{C}_{\chi_\bomega}$.
Notice that either $$V(\bomega_0)\subset \bu v\ \  or \ \ \bu
v\cap V(\bomega_0)=0.$$  If $\chi_\bomega\ne\chi_{\bomega_0}$ then
the second possibility must hold and so  $$V\cong
V(\bomega_0)\oplus \bu v.$$ This shows that induction begins.

Now assume that $U$ is indecomposable with length $\ell> 1$ and that we
know the result for all modules with length strictly smaller than $\ell$.
Let $U_1$ be a proper non--trivial
submodule of $U$ and consider the short exact sequence,
\begin{equation*}
0\to U_1\to U\to U_2 \to 0
\end{equation*}
Since  $\text{Ext}^1_{\cal C_q}(U_j,V(\bomega_0)) = 0$ for $j=1,2$ by
the induction hypothesis, the result follows by using the exact
sequence $\text{Ext}^1_{\cal C_q}(U_2,V(\bomega_0))\to
\text{Ext}^1_{\cal C_q}(U,V(\bomega_0))\to
\text{Ext}^1_{\cal C_q}(U_1,V(\bomega_0))$. Part (ii) is now
immediate by using a similar induction on the length of $V_2$.
\end{pf}

\subl{}
The proof of Theorem \ref{blockdec}(i) is completed as follows. Let
$V$ be an indecomposable $\bu$-module. We prove that there exists
$\chi\in\Xi_q$ such that $V\in\cal C_\chi$  by an induction on the
length of $V$. If $V=V(\bomega)$ is irreducible it follows from
the definition of $\cal{C}_{\chi_\bomega}$.  If $V$ is reducible,
let $V(\bomega_0)$ be an irreducible subrepresentation of $V$ and
let $U$ be the corresponding quotient. In other words, we have an
extension
\begin{equation*}
0\to V(\bomega_0)\to V\to U \to 0
\end{equation*}
Write $U=\oplus_{j=1}^r U_j$ where each $U_j$ is indecomposable.
By the inductive hypothesis, there exist $\chi_j\in \Xi_q$ such that
$U_j\in\cal C_{\chi_j}$, $1\le j\le r$. Suppose that there exists
$j_0$ such that $\chi_{j_0}\ne \chi_{\bomega_0}$. Lemma
\ref{nullext} implies that $$\text{Ext}^1_{\cal C_q}(U,V(\bomega_0))
\cong \oplus_{j=1}^r \rm{Ext}^1_{\cal C_q}(U_j,V(\bomega_0))\cong
\oplus_{j\neq j_0} \rm{Ext}^1_{\cal C_q}(U_j,V(\bomega_0)).$$ In other
words, the exact sequence $0\to V(\bomega_0)\to V\to U \to 0$ is
equivalent to one of the form
\begin{equation*}
0\to V(\bomega_0)\to  V'\oplus U_{j_0}\to \oplus_{j\neq
j_0} U_j \oplus U_{j_0}\to 0
\end{equation*}
where
\begin{equation*}
0\to V(\bomega_0)\to V'\to \oplus_{j\neq j_0} U_j \to 0
\end{equation*}
is an element of $\oplus_{j\neq j_0}
\rm{Ext}^1_{\cal C_q}(U_j,V(\bomega_0))$. But this contradicts the fact
that $V$ is indecomposable. Hence $\chi_j=\chi_{\bomega_0}$ for
all $1\le j\le r$ and $V\in\cal C_{\chi_{\bomega_0}}$.
\hfill\qedsymbol

\subl{} We now prove Theorem \ref{blockdec}(ii). The idea is
similar to the one used  in \cite{EM}, although, again, with
the theory of $\ell$-lattices the proof is simpler,  uniform and
works for generic $q$. The proof we give depends on Proposition
\ref{gammaq} which has only been stated so far for the classical
Lie algebras. The proof of part (ii) of Theorem \ref{blockdec} for the
exceptional Lie algebras is postponed to the appendix after we
state the analog of Proposition \ref{gammaq} for these algebras.

We first consider the cases $\lie g=A_n,B_n,C_n,D_n$, where,
 in the case of $D_n$, we assume that $n$ is odd. Thus let $i_\bullet$ be the unique element in
$I_\bullet$.  From now on we denote by $V(a)$ the irreducible
fundamental representation $V(\bomega_{i_\bullet,a})$.

We shall need the following result.
{\samepage
\begin{prop}\label{catgen}\hfill
\begin{enumerit}

 \item[(i)]   Let $a_1,\cdots ,a_k\in\bc^\times$ and let $\sigma$ be
any permutation of $\{1,\cdots ,k\}$. Then, the modules
$V(a_1)\otimes\cdots\otimes V(a_k)$ and
$V(a_{\sigma(1)})\otimes\cdots\otimes V(a_{\sigma(k)})$ are
linked.
\item[(ii)]Given $\bomega\in\cal{P}_q^+$, there exists a set
(possibly not unique)
   $S_\bomega=\{a_1,\cdots ,a_k\}\subset\bc^\times
$ such that
 $W(\bomega)$ is a subquotient  of $V(a_1)\otimes\cdots\otimes
 V(a_k)$ and hence $W(\bomega)$ and $V(\bomega)$ are  linked to it.
\end{enumerit}\hfill\qedsymbol
\end{prop}}
\begin{pf}    To prove (i), note first that since the Grothendieck ring
of the category of finite--dimensional representations is
commutative, \cite{FR}, it follows that the   modules
$V(a_1)\otimes\cdots\otimes V(a_k)$ and
$V(a_{\sigma(1)})\otimes\cdots\otimes V(a_{\sigma(k)})$ have the
same irreducible constituents for all permutations $\sigma$ of
$\{1,2,\cdots ,k\}$. Hence to show that they are linked it
suffices to prove that there exists a permutation $\tau$ such that
$V(a_{\tau(1)})\otimes\cdots\otimes V(a_{\tau(k)})$ is
indecomposable. But this is clear using Theorem \ref{cyc}, which
implies that there  exists a
 permutation  $\tau$ of $a_1,\cdots ,a_k\in\bc^\times$ such that
 $V(a_{\tau(1)})\otimes\cdots\otimes V(a_{\tau(k)})$ is cyclic
 on the tensor product of highest
 weight vectors and hence indecomposable. By Theorem
\ref{cpconj}   it suffices to prove (ii) when
$\bomega=\bomega_{i,a}$ for some $i\in I$, $a\in \bc^\times $. But
 this follows from
 \cite[Theorem 6.1, Proposition 7.5, Theorem 8.2]{CPdorey}. It is
 immediate from (i) that $W(\bomega)$ is linked to
 $V(a_1)\otimes\cdots\otimes V(a_k)$.
\end{pf}

\begin{cor} Let $\bomega,\bomega'\in\cal P_q^+$ and assume that
$S_\bomega=\{a_1,\cdots ,s_k\}$, $S_{\bomega'}=\{a_1', \cdots
,a_l'\}$. Then $V(\bomega)$ and $V(\bomega')$ are linked iff
$V(a_1)\otimes\cdots\otimes V(a_k)$ and
$V(a_1')\otimes\cdots\otimes V(a_l')$ are linked.
\end{cor}

As a consequence to prove
 Theorem \ref{blockdec}(ii) it suffices to show that, if the modules
  $V(a_1)\otimes\cdots\otimes V(a_k)$ and $V(b_1)\otimes\cdots\otimes V(b_s)$
   have the same elliptic character, then they are linked.


\subl{} The next result identifies minimal sets $S_\boe$, where
$\boe$ is the identity element in $\cal{P}_q^+$. Notice that in
this case the associated irreducible representation is the trivial
one.

\begin{prop}\label{minsets} For all $a\in\bc^\times$ we can take,
\begin{alignat*}{3}
S_\boe(a) & = \{a, aq^{2}\cdots ,aq^{2n}\},& \text{if } \lie g = A_n, \\
& = \{a,aq^{4n-2}\},& \text{ if } \lie g= B_n,\\
& = \{a,aq^{2n+2}\}, & \text{ if } \lie g = C_n,\\
& = \{a,aq^2,aq^{2n-2},aq^{2n}\}, & \text{ if } \lie g = D_{n}.
\end{alignat*}
\end{prop}

\begin{pf} It was proved in \cite[Proposition 5.1]{CPminaff1}
that the dual  $V(\bomega)^*$ is isomorphic to $V(\bomega^*)$
where $(\bomega^*)_i=(\bomega)_{w_0i}(q_i^{\check h}u)$. This
proves the statements for $B_n$ and $C_n$. For $A_n$, it suffices
to prove that $$V(\bomega_{i,q^{i-1}})\subset
V(\bomega_{1,q^i})\otimes V(\bomega_{i-1,1}),\ \ 1\le i\le n.$$
But this follows from \cite[Lemma 5.2]{CPminaff1}. The proof for
$D_n$ is similar and we omit the details.
\end{pf}

\subl{} Given $a\in\bc^\times$, let $\bomega_a\in\cal{P}_q^+$ be the element
defined by
$$(\bomega_a)_j=1,\ \ j\ne i_\bullet,\ \
(\bomega_a)_{i_\bullet}=\prod_{a_j\in S_{\boe}(a)}(1-a_ju).$$
The following is now immediate.
\begin{cor} Let $V'\in\cal{C}_q$, then $V'$ is linked to
$V'\otimes
W(\bomega_a).$
\end{cor}

\subl{}\label{pbd}  Recall from Section \ref{gammaq} that the group $\Xi_q$ is isomorphic to the quotient of the free group $\tilde\Xi_q$,
generated  by  elements $\chi_a$, by the subgroup generated by
 $$\kappa_a=\sum_{a_j\in
S_\boe(a)}\chi_{a_j}$$
for all $a\in\bc^{\times}$. We can now complete the proof of Theorem
\ref{blockdec}(ii). Suppose that the modules
$V(a_1)\otimes\cdots\otimes V(a_k)$ and
$V(b_1)\otimes\cdots\otimes V(b_s)$ have the same elliptic
character.
 Then, in $\tilde\Xi_q$, we have
$$\sum_{r=1}^k \chi_{a_r} - \sum_{r=1}^s \chi_{b_r} = \sum_{r=1}^p
m_r \kappa_{c_r}$$ for some $c_r\in \bc^{\times}$ and integers
$m_r$. We can assume that there exists $p'$ such that $m_r\le 0$ if $1\le r\le p'$ and
$m_r>0$ otherwise. Now we have, $$\sum_{r=1}^k
\chi_{a_r}+\sum_{r=1}^{p'}(-m_r) \kappa_{c_r} = \sum_{r=1}^s
\chi_{b_t} + \sum_{r=p'+1}^p m_r \kappa_{c_r}.$$ Since this is an
equality in a free group, it follows that we have an equality of
sets with multiplicities,
$$\{a_1,\cdots ,a_k\} \cup\cup_{r=1}^{p'} S_{\boe}(c_r)=\{b_1,\cdots
,b_s\}\cup\cup_{r=p'+1}^{p} S_{\boe}(c_r),$$
where the multiplicity of $S_{\boe}(c_r)$ is $-m_r$ if $r\leq p'$, and $m_r$ if $r>p'$.

This now gives,
\begin{align*}
V(a_1)\otimes\cdots\otimes V(a_k) \sim &\ V(a_1)\otimes\cdots\otimes V(a_k)
\otimes(\otimes_{r=1}^{p'}W(\bomega_{\boc_r})^{\otimes (-m_r)})\\
\sim\ & V(b_1)\otimes\cdots\otimes
V(b_s)\otimes(\otimes_{r=p'+1}^p W(\bomega_{\boc_r})^{\otimes m_r})\\
\sim\ & V(b_1)\otimes\cdots\otimes V(b_s)
\end{align*}
where $\sim$ stands for linking relation.
The proof of Theorem \ref{blockdec}(iii) is immediate.

\subl{} We now consider the $\lie g =D_n$ for even $n$.
 Thus, set $V_-(a)=V(\bomega_{n-1,a})$ and $V_+(a)=V(\bomega_{n,a})$. We state the analogue of Proposition \ref{catgen} which is proved in a similar way.

\begin{prop}\hfill
\begin{enumerit}
 \item[(i)]  Let $a_1,\cdots ,a_{k}\in \bc^\times, \varepsilon_1,\cdots,\varepsilon_{k}\in \{+,-\}$, and $\sigma$ be a permutation of $\{1,\cdots,k\}$. Then the modules $V_{\varepsilon_1}(a_1)\otimes\cdots \otimes V_{\varepsilon_k}(a_{k})$ and $V_{\varepsilon_{\sigma(1)}}(a_{\sigma(1)})\otimes\cdots \otimes V_{\varepsilon_{\sigma(k)}}(a_{\sigma(k)})$ are linked.

\item[(ii)]  Given $\bomega\in\cal{P}_q^+$, there exists a (non--unique) pair of sets
$S_\bomega=(\{a_1,\cdots ,a_{k}\},\{\varepsilon_1,\cdots,\varepsilon_{k}\}),$
where $a_j\in \bc^\times$ and $\varepsilon_j\in \{+,-\}$, such that $W(\bomega)$ is a subquotient of $U=V_{\varepsilon_1}(a_1)\otimes\cdots \otimes V_{\varepsilon_k}(a_{k})$. In particular $W(\bomega)$ and $V(\bomega)$ are linked to $U$.
\end{enumerit}\hfill\qedsymbol
\end{prop}

The analogue of Corollary \ref{catgen} is immediate.

\subl{} Given $a\in\bc^{\times}$, define the sets

\begin{align*}
S_{\boe,0}(a) = \{a, aq^{2n-2}\}, \qquad S_{\boe,1}(a)  =\{a,aq^2\}
\end{align*}
and, for $k\in\{0,+,-\}$, let $\bomega_{a,k}\in\cal P_q$ be given by
$$(\bomega_{a,k})_j=1,\ \ j\notin I_\bullet$$
and
$$(\bomega_{a,k})_{n-1}=
\begin{cases}
1, &\text{ if } k=+,\\
\prod\limits_{a_j\in S_{\boe,0}(a)}(1-a_ju), &\text{ if } k=-,\\
\prod\limits_{a_j\in S_{\boe,1}(aq^{2n-2})}(1-a_ju), &\text{ if } k=0,
\end{cases}
\quad\quad (\bomega_{a,k})_{n}=
\begin{cases}
1, &\text{ if } k=-,\\
\prod\limits_{a_j\in S_{\boe,0}(a)}(1-a_ju), &\text{ if } k=+,\\
\prod\limits_{a_j\in S_{\boe,1}(a)}(1-a_ju), &\text{ if } k=0.
\end{cases}$$

The next Proposition is proved exactly like Proposition \ref{minsets}.

\begin{prop}\label{minsetsDneven}
The trivial representation is linked to $W(\bomega_{a,k})$.
\end{prop}

\begin{cor} Let $V'\in\cal{C}_q$, then $V'$ is linked to
$V'\otimes W(\bomega_{a,k}).$
\end{cor}

\subl{} We complete the proof of Theorem \ref{blockdec}(ii) as  in
section \ref{pbd}  using Proposition \ref{gammaq}.

\section{Appendix: Exceptional Algebras}

We now consider the problem of determining the block decomposition when $\lie g$ is one of the exceptional algebras.

\subl{} Thus, let $i=1$ and $V=V(\bomega_{1,a})$.
We start observing that Proposition and Corollary \ref{catgen} holds for the exceptional algebras.

\subl{} Given $a\in \bc^{\times}$, define sets $S_{\boe,k}(a)$, where $k=1,2,$ when $E_6,E_7, F_4, G_2$,  $k=1,2,3$, when $\lie g=E_8$, as follows.
\begin{alignat*}{3}
&S_{\boe,1}(a)  =\{a,aq^8,aq^{16}\}, \ S_{\boe,2}(a)  = \{a,aq^2,aq^4,aq^{12},aq^{14},aq^{16}\}, & &\ \text{if } \lie g = E_6,\\
&S_{\boe,1}(a)  =\{a,aq^{18}\}, \ S_{\boe,2}(a)  = \{a,aq^2,aq^{12},aq^{14},aq^{24},aq^{26}\}, & &\ \text{if } \lie g = E_7,\\
&S_{\boe,1}(a)  =\{a,aq^{30}\}, \ S_{\boe,2} = \{a,aq^{20},aq^{30}\},\ S_{\boe,3}(a)  = \{a,aq^{12},aq^{24},aq^{36},aq^{48}\}, & &\ \text{if } \lie g = E_8,\\
&S_{\boe,1}(a)  =\{a,aq^{18}\}, \ S_{\boe,2}(a)  = \{a,aq^2,aq^{12},aq^{24}\}, & &\ \text{if } \lie g = F_4,\\
&S_{\boe,1}(a)  =\{a,aq^{12}\}, \ S_{\boe,2}(a)  = \{a,aq^8,aq^{16}\}, & &\ \text{if } \lie g = G_2,\\
\end{alignat*}

Now, define $\bomega_{a,k}\in \cal P_q^+$ by
$$(\bomega_{a,k})_j=1,\ \ j\ne 1,\ \ (\bomega_{a,k})_1=\prod_{a_j\in S_{\boe,k}(a)}(1-a_ju).$$

Similarly to Proposition \ref{minsets}, we have the following.

\begin{prop}
The trivial representation is linked to $W(\bomega_{a,k})$.
\end{prop}

\begin{cor} Let $V'\in\cal{C}_q$, then $V'$ is linked to
$V'\otimes W(\bomega_{a,k}).$
\end{cor}

\subsection{}  Now we complete the statement of Proposition \ref{gammaq} for the exceptional algebras.
\begin{prop}\label{gammaqe}
Assume that $\lie g$ is of type $E_6,E_7,E_8,F_4$ or
$G_2$.  The group $\cal{P}_q/\cal{Q}_q$ is isomorphic to the
(additive) abelian group $\Xi_q$ with generators
$\{\chi_a:a\in\bc^\times\}$ and relations:
$$\sum_{a_j\in S_{\boe,k}(a)}\chi_{a_j}=0$$
for all $a\in\bc^\times$ and $k=1,2$, for $\lie g\ne E_8$ and $k=1,2,3$, for $\lie g=E_8$.
\end{prop}

The idea of the proof is the same of Proposition \ref{gammaq} and requires a long, but straight forward, case by case checking.

\subl{} The rest of the proof of Theorem \ref{blockdec} (ii) and
(iii) for the exceptional algebras is similar to section \ref{pbd}
using Proposition \ref{gammaqe}.

\bibliographystyle{amsplain}

\end{document}